\pdfoutput=1
\documentclass[]{article}
\usepackage[a4paper,margin=25mm]{geometry}
\usepackage{graphicx}
\usepackage{amsmath}
\usepackage{amssymb}
\usepackage{hyperref}
\usepackage{listings}  
\usepackage{tabu,multirow}
\usepackage{amsfonts}

\usepackage{pgfplots}
\pgfplotsset{compat=1.15}
\PassOptionsToPackage{pdf}{pstricks}
\usepackage{pstricks,pst-func} 
\usepackage{tikz}

\usepackage{mathrsfs}
\usetikzlibrary{arrows}
\usepackage[thinc]{esdiff}

\DeclareMathSymbol{\shortminus}{\mathbin}{AMSa}{"39}

\usetikzlibrary{decorations.pathreplacing,calc}
\newcommand{\tikzmark}[1]{\tikz[overlay,remember picture] \node (#1) {};}

\usepackage{algpseudocode}
\usepackage{algorithm}

\usepackage{fancyhdr}
\title{Rational right triangles and the Congruent Number Problem}
\author{G. Jacob Martens}
\begin{document}
\definecolor{uuuuuu}{rgb}{0.26666666666666666,0.26666666666666666,0.26666666666666666}
\definecolor{xdxdff}{rgb}{0.49019607843137253,0.49019607843137253,1}
\definecolor{ududff}{rgb}{0.30196078431372547,0.30196078431372547,1}
\definecolor{qqccqq}{rgb}{0,0.8,0}
\definecolor{ffqqtt}{rgb}{1,0,0.2}
\definecolor{qqzzcc}{rgb}{0,0.6,0.8}
\algnewcommand{\algorithmicgoto}{\textbf{go to}}%
\algnewcommand{\Goto}[1]{\algorithmicgoto~\ref{#1}}%

\maketitle

\small

\begin{abstract}
From Euclid's fundamental formula for the Pythagorean triples we define the rational triples relating certain congruent numbers by an identity and explore their relationships. We introduce two geometric methods relating the congruent number problem to pairs of conic sections. We show a relationship between the Cassini ovals and the congruent number problem. By the tangent method we define a set of rational triangles from an initial solution for a congruent number. We define the prime footprint equations for right triangles for certain congruent numbers. By the unseen recurrence we define infinite trees of rational triangles and congruent numbers. Congruent number families are defined related to the Fibonacci/Lucas numbers and the Chebyshev polynomials. We show that the semiperimeter of the Brahmagupta triangles are congruent numbers in function of the Chebyshev polynomials of the first kind. We present a naive recursive algorithm to solve the problem posed by Fermat for triangles with $a+b=\square$ and $c=\square$
\end{abstract}
\mbox{}
\newline
\linebreak
\mbox{}
\newline
\linebreak
\mbox{}
\newline
\linebreak
\includegraphics[width=155mm,scale=1]{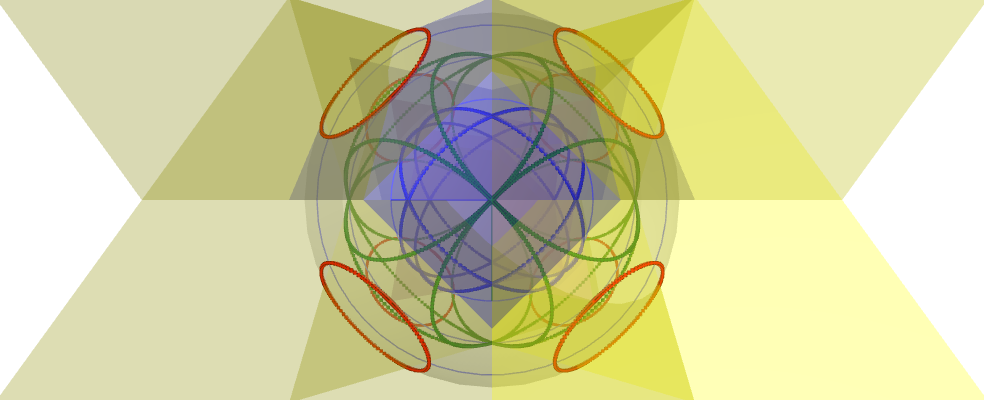}
\newline
\linebreak
\begin{center}
For 
\linebreak
\newline
Dr. N.Tati Ruiz
\linebreak
\newline
A moment of opportunity...
\end{center}
\mbox{}
\newline
\linebreak
\begin{center}
This pre-publication version is free to view and download for personal use only. 
\end{center}
\begin{center}
Not for re-distribution, re-sale or use in derivative works. 
\end{center}
\begin{center}
\copyright $\;$ G. Jacob Martens 2021.
\end{center}
\pagebreak
\tableofcontents
\pagebreak
\section{Preface}
This mathematical journey started by exploring a problem formulated by Fermat \cite{1} to find a right triangle having integer sides $a,b,c$ such that $a+b=\square$ and $c=\square$.
In $1643$ Fermat wrote a letter to Mersenne confirming the smallest solution to be (a,b,c)=(456548602761,1061652293520,4687298610289).
To obtain this result he must have had a rigorous method, first to find this solution and second to claim that it is the smallest triangle of this kind. 
My own survey and solution to this problem is outlined in the last section. 
\newline
\linebreak
It led me to see the importance of the fundamental theorem of Pythagoras which one can compare with a swiss army knife still hiding multiple relationships and identities.
The title was an obvious choise since most of the discoveries relate to the Congruent Number Problem \cite{2}, one of the oldest unsolved problems in mathematics, which is about trying to find a right triangle with rational sides whose area is a certain positive integer\cite{3}.
\newline
\linebreak
Starting from Euclid's fundamental formula for the Pythagorean triples \cite{4} we define the rational triples from which we will derive several identities and relationships.
The first identity defines a relationship between the areas of the rational triples which are all congruent numbers.
A second identity shows how the sum of the squares of the differences of two sides relate to the Euclidean distance \cite{5} between two points in a three-dimensional Euclidean space. 
Another relationship defines three sets of quadratic Diophantine equations also known as Eulers Concordant Forms \cite{6} for the congruent number problem.
By another identity, relating the sum of squares of the sides, we define the Trinity system, a three-dimensional system of spheres embedding small and great circles,  preserving the Pythagorean equation accross the spheres. In this context we define the Trinity vectors, their derivatives and their product relationships defining parallelepipeds with volume equal to $1/2$.
\newline
\linebreak
A geometric method shows how congruent numbers are related to two conjugate conics and proves that the rational points on the congruent number elliptic curve $E_N: y^2=x^3-N^2x$ are related to rational intersection points of a line and an ellipse.
In this context we show how the famous often celebrated solution found  by Don Zagier for the congruent number $157$  can be represented by the two numbers $87005 ,610961$.
We define and apply the line ellipse intersection method to prove that the polynomial $(4t^2+1)(4t^2-8t+5)$ represents congruent numbers.
Using this method we also define congruent number polynomials from lattice points of an ellipse.
\newline
\linebreak
We define the twin hyperbolas which is another method relating two conic sections to the congruent number problem and show how their rational intersection points define the two congruent number polynomials
 $ 2 \left(11 t^4-36 t^3+30 t^2-12  t+19\right)$ and $2\left(11 t^4+60 t^3+66 t^2-132  t+43\right)$.
\newline
\linebreak
We define a relationship between the Cassini oval \cite{7} axis intersection points and a system of equations defined by Kurt Heegner in 1952 \cite{8} related to the congruent number problem and show how to obtain Cassini ovals with two or four $X$-axis intersections for the same congruent number.
\newline
\linebreak
A method named the tangent method, due to its relationship to the point doubling method on elliptic curves, shows how a set of rational right triangle solutions can be generated from the known rational hypotenuse of an initial solution for one congruent number by solving a binary quadratic form.
\newline
\linebreak
We define the prime footprint equations which are a template for prime or almost prime congruent numbers and include tables with their solutions for numbers less than $1000$. 
From these emperically defined equations one can derive relationships to conic sections and their rational points.
\newline
\linebreak
We show how a hidden recurrence defines infinite trees of congruent number sequences from an initial Pythagorean triple.
By defining a method choosing the triple tree-side-walk it is possible to generate different congruent number trees for the chosen side-path.
Explicit examples for the side-paths $a^i,\;ab,\;bb$ are given for Euclid's definition.
\newline
\linebreak
A method is given to define congruent number families related to the Fibonacci\cite{9}/Lucas\cite{10} numbers and the Chebyshev polynomials \cite{11}. 
For both we look at some properties of the associated congruent number elliptic curves defining some general rational points.
We show how the solution for a right triangle with area $5$ first found by Fibonacci is related to the Fibonacci numbers. 
We also show that the Brahmagupta triangles \cite{12} can be defined in function of the Chebyshev polynomials and that their semiperimeters are congruent numbers.\newline
For these triangles we also define a general elliptic curve and review some of its basic properties. 
\newline
\linebreak
We conclude with the naive recursive algorithm generating the solutions for the Fermat triangles which ignited this mathematical journey.
\newline
\linebreak
In order to verify most results the following Sagemath \cite{13} notebooks are published along with this paper.
\begin{align*}
 &\texttt{02\_Beyond\_Pythagoras.ipynb} \\
 &\texttt{03\_Conjugate\_conics.ipynb} \\
 &\texttt{04\_Twin\_hyperbolas.ipynb} \\
 &\texttt{05\_Cassini\_ovals.ipynb} \\
 &\texttt{06\_Tangent\_method.ipynb} \\
 &\texttt{07\_Prime\_footprints.ipynb} \\
 &\texttt{08\_Unseen\_recurrence.ipynb} \\
 &\texttt{09\_Fibonacci\_Lucas.ipynb} \\
 &\texttt{10\_Chebyshev\_polynomials.ipynb} \\
 &\texttt{11\_Fermat.ipynb}
\end{align*}
\newline
\linebreak
\includegraphics[width=155mm]{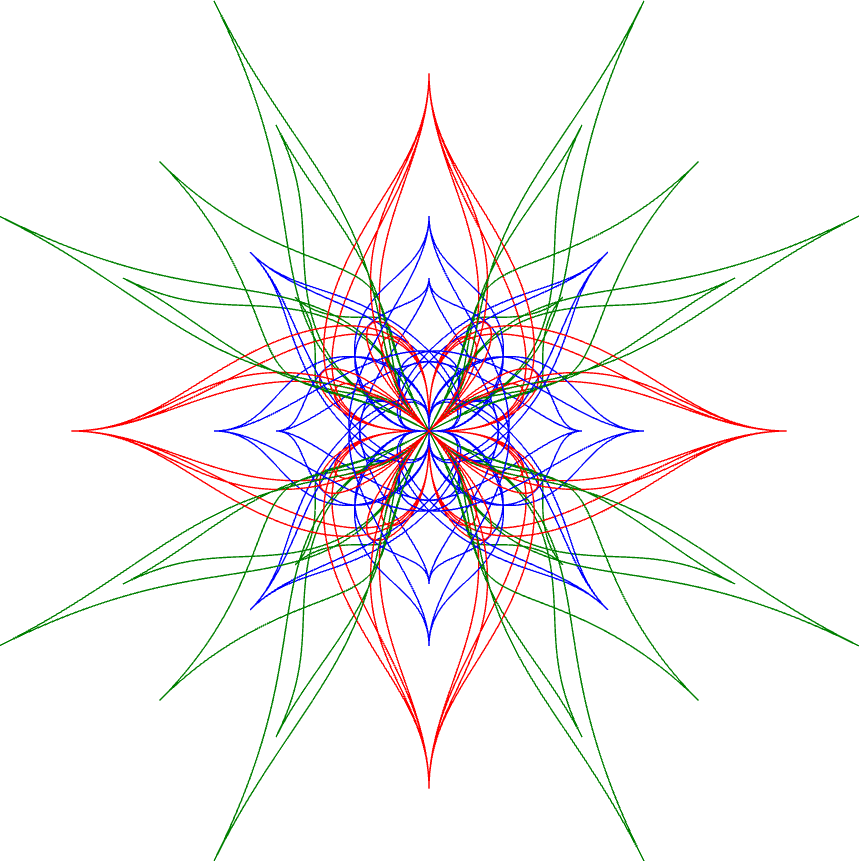}
\pagebreak
\section{Beyond Pythagoras }
From Euclid's fundamental formula we define the rational triples for right angled triangles.
\newline
The area's of these right triangles are congruent numbers between which certain relationships and identities exist.
\newline
The rational triples are also the foundation for the Trinity system a three-dimensional differential system defining rational points on circles of spheres and vector algebraic relationships. 
\subsection{The rational triples }
The Pythagorean triples $A,B,C$  representing a right triangle with integer sides and area $N=\frac{AB}{2}$ , satisfying the equation $A^2+B^2=C^2$  , are generated by Euclid's fundamental formula 
\begin{equation} (A,B,C)=\left(m^2-n^2, 2 m n, m^2+n^2\right) \;\; \forall \;\; { \; m,n \in  \mathbb{Z} \; and \; m>n>0 }\end{equation}
Combining two sides we can define three related right triangles satisfying $a'^2+b'^2= c'^2$ 
\small
\begin{eqnarray}  \nonumber
(a, b, c)_{AC}' & = & \left( 2 A^2 C^2, B^2 (A^2 + C^2), A^4 + C^4 \right)\\  \nonumber
(a, b, c)_{BC}' &=& \left( 2 B^2 C^2, A^2 (B^2 + C^2), B^4 + C^4\right) \\
(a, b, c)_{BA}' &=& \left( 2 B^2 A^2, C^2 (B^2 - A^2), B^4 + A^4\right)
\end{eqnarray} 
Reducing the sides by $ABC$ we obtain the rational triples $a,b,c$ 
 satisfying the equation $a^2+b^2=c^2\;\; \forall \;\; a,b,c \; \in  \mathbb{Q}$ 
\small
\begin{eqnarray}   \nonumber
(a,b,c)_{AC} =\left(\frac{2 A C}{B},\frac{B (A^2 +C^2)}{A C},\frac{A^4+C^4}{A B C}\right)
=&\left( \frac{m^4-n^4}{mn},\frac{4mn(m^4+n^4)}{m^4-n^4},\frac{4m^4n^4+\left(m^4+n^4\right)^2}{nm(m^4-n^4)}\right)  \\  \nonumber
  (a,b,c)_{BC}  =\left(\frac{2 B C}{A},\frac{A (B^2 +C^2)}{B C},\frac{B^4+C^4}{A B C}\right)   
 =&\left(\frac{4mn(m^2+n^2)}{m^2-n^2},\frac{(m^2-n^2)\left(4m^2n^2+\left(m^2+n^2\right)^2\right)}{2 m n (m^2+n^2)},\frac{16 m^4 n^4 + (m^2+n^2)^4}{2mn(m^4-n^4)}\right)  \\
 (a,b,c)_{BA}  =\left(\frac{2 B A}{C},\frac{C (B^2 -A^2)}{B A},\frac{B^4+A^4}{A B C}\right)   
 = &\left(\frac{4mn(m^2-n^2)}{m^2+n^2},\frac{(m^2+n^2)\left(4m^2n^2-\left(m^2-n^2\right)^2\right)}{2 m n (m^2-n^2)},\frac{16 m^4 n^4 + (m^2-n^2)^4}{2mn(m^4-n^4)}\right) 
 \end{eqnarray}
 \subsection{The area identity}
The area's for the four triangles $(A,B,C),(a,b,c)_{AC,BC,BA}$  represented by 
\small
\begin{eqnarray}  \nonumber
(N,N_{AC},N_{BC},N_{BA})&=& \left(\frac{AB}{2},A^2+C^2,B^2+C^2,B^2-A^2\right)  \\
&= &\left(m n \left(m^2-n^2\right),2 \left(m^4+n^4\right),m^4+6 m^2 n^2+n^4,-m^4+6 n^2 m^2-n^4\right) 
\end{eqnarray}
are congruent numbers related by the identity
\begin{equation}   N_{AC}^2 + N_{BC}^2 + N_{BA}^2= 6 \left(C^4 - 4 N^2\right)  \end{equation}
For example for $(m,n)=(2,1)$, we get the right triangles 
\small
$$ (A,\;B,\;C)=(3, 4, 5) \;\;,\;\;
 (a,b,c)_{AC}=\left(\frac{15}{2}, \frac{136}{15}, \frac{353}{30}\right) ,\;
 (a,b,c)_{BC}=\left(\frac{40}{3}, \frac{123}{20}, \frac{881}{60}\right) ,\;
 (a,b,c)_{BA}=\left(\frac{24}{5}, \frac{35}{12}, \frac{337}{60}\right)$$
relating the congruent numbers     
$ (N,N_{AC},N_{BC},N_{BA})=(6,34,41,7)$
as follows 
$$ 34^2+41^2+7^2=6 (5^4- 4*6^2)$$
\subsection{The connecting line}
The parallel lines to the $x$-axis at $y$ values $\pm2ABC$ intersect and connect the congruent number elliptic curves 
\begin{equation} 
 E_{AC}: y^2 = \; x^3 -N_{AC}^2\;x  \;,\; E_{BC}: y^2 = \; x^3 -N_{BC}^2\;x  \;,\;E_{BA}: y^2 = \; x^3 -N_{BA}^2\;x  
 \end{equation}
 at the rational points $P_{AC,BC,BA}:(x,y)$
\begin{align}   \nonumber
 P_{AC} : &\left(-B^2,\pm2ABC\right)  =  \left(\shortminus4 m^2 n^2 \;\;\;\;\;\;\;,  \pm 4 m n (m^2 - n^2) (m^2 + n^2) \right) \\  \nonumber
 P_{BC} : & \left(-A^2,\pm2ABC\right) =   \left(\shortminus(m^2 - n^2)^2,  \pm 4 m n (m^2 - n^2) (m^2 + n^2) \right) \\
 P_{BA} : & \left(\;\;\;C^2,\pm2ABC\right)  =  \left( (m^2 + n^2)^2\;,  \pm4 m n (m^2 - n^2) (m^2 + n^2) \right) 
\end{align}
For the example $(m,n)=(2,1)$, we get the three positive and negative colinear rational points 
$$ P_{AC} :(-16,\pm120) \;\;\;\;, \;\;\;\;\;  P_{BC} : (-9,\pm120) \;\;\;\;, \;\;\;\;\; P_{BA} : (25,\pm120) $$
\subsection{The Diophantine equations} 
Each Pythagorean triple $(a,b,c)_{AC,BC,BA}$  forms a pair of quadratic Diophantine equations, also known as Euler's concordant forms for the congruent number problem
\begin{equation}x^2+N y^2 =z^2 \; and  \;x^2-N y^2 = t^2 \;\; \forall \;\; { \; x,y,z,t,N} \in  \mathbb{Z}\end{equation}
relating the  numerator and the denominator of the hypotenuse $c$  to the area $N_{AC,BC,BA}$ such that  $$(x,y)= (Numerator[c],2*Denominator[c])$$
with the following parameterization 
\small
\begin{eqnarray} 
(x,y)_{AC}  &=& (A^4 + C^4 ,2 A B C ) = \left(2\left (m^8 + 6 m^4 n^4 + n^8\right), 4 nm(m^4-n^4) \right) \\ \nonumber
(z,\;t)_{AC} & = & \left(\sqrt{(A^4 + C^4)^2 + (A^2 + C^2) (2 A B C )^2}  ,\sqrt{(A^4 + C^4)^2 - (A^2 + C^2) (2 A B C )^2} \right) \\  
&= & \left(2(m^8 + 4 m^6 n^2 - 2 m^4 n^4 + 4 m^2 n^6 + n^8),2(m^8 - 4 m^6 n^2 - 2 m^4 n^4 - 4 m^2 n^6 + n^8)\right) \\[0.5cm]
(x,y)_{BC}  &=& (B^4 + C^4 ,2 A B C ) = \left(m^8 + 4 m^6 n^2 + 22 m^4 n^4 + 4 m^2 n^6 + n^8, 4 nm(m^4-n^4) \right) \\ \nonumber
(z,\;t)_{BC} & = & \left(\sqrt{(B^4 + C^4)^2 + (B^2 + C^2) (2 A B C )^2}  , \sqrt{(B^4 + C^4)^2 - (B^2 + C^2) (2 A B C )^2} \right) \\  
&=& \left(m^8 + 12 m^6 n^2 + 6 m^4 n^4 + 12 m^2 n^6 + n^8,m^8 - 4 m^6 n^2 - 26 m^4 n^4 - 4 m^2 n^6 + n^8\right) \\[0.5cm]
(x,y)_{BA}  &=& (B^4 + A^4 ,2 A B C )= \left(m^8 - 4 m^6 n^2 + 22 m^4 n^4 - 4 m^2 n^6 + n^8, 4 nm(m^4-n^4) \right) \\ \nonumber
(z,\;t)_{BA} & = &  \left(\sqrt{(B^4 + A^4)^2 + (B^2 - A^2) (2 A B C )^2}  ,\sqrt{(B^4 + A^4)^2 - (B^2 - A^2) (2 A B C )^2} \right) \\  
& = & \left(m^8 - 12 m^6 n^2 + 6 m^4 n^4 - 12 m^2 n^6 + n^8,m^8 + 4 m^6 n^2 - 26m^4 n^4 + 4 m^2 n^6 + n^8\right)
\end{eqnarray} 
This results in three pairs of equations for which the $y$ solutions are all equal to $2ABC$
\begin{equation}\left(y^2 = \frac{z^2 - x^2}{N} = \frac{x^2 - t^2}{N}\right)_{AC,BC,BA} \;\;\; and \;\;\; y_{AC}= y_{BC}= y_{BA}=2 A B C\end{equation}
For the example $(m,n)=(2,1)$ we get the following solutions 
$$(x,y,z,t,N)_{AC}=(706, 120, 994, 94,34)$$
$$\;\;\;(x,y,z,t,N)_{BC}=(881, 120, 1169, 431, 41)$$
$$(x,y,z,t,N)_{BA}=(337, 120, 463, 113, 7)$$
And verifying the identity $(15)$  we find for $y_{AC,BC,BA}=120$ that  
$$120^2 = \frac{994^2 - 706^2}{34}= \frac{706^2 - 94^2}{34} = \frac{1169^2 - 881^2}{41}= \frac{881^2 - 431^2}{41} = \frac{463^2 - 337^2}{7}= \frac{ 337^2 - 113^2}{7}$$
 \subsection{The distance identity}
A distance identity closely related to the Euclidean distance formula for three-dimensions, also known as the 2-Norm is obtained by enumerating the triples $(a,b,c)_{AB,BC,BA}$ as cartesian coordinates $(a,b,c)_{1,2,3}$.
\newline
Taking the sum of the squares of the difference of the sides $a,c$ we obtain the identity
\begin{eqnarray}\left(2 \frac{C^4 - 3 (A B)^2 }{A B C}\right)^2 
& = & 2 \left(\left(c_1 - a_1\right)^2 + \left(c_2 - a_2\right)^2 + \left(c_3 - a_3\right)^2\right) \\
& = &\;\;\;\left(\left(c_1 - a_1\right) \;+ \left(c_2 - a_2\right) \;\;+ \left(c_3 - a_3\right) \right)^2 \\
& = & 4\left(\left(c_1 - a_1\right) \left(c_2 - a_2\right) + \left(c_1 - a_1\right) \left(c_3 - a_3\right)+\left(c_2 - a_2\right) \left(c_3 - a_3\right) \right)
\end{eqnarray}
for which the three products of the two differences of $(18)$ are also perfect rational squares 
\begin{align}  \nonumber
4 (c_1 - a_1) (c_2 - a_2) &=  \;\left(\;\left(c_1 - a_1\right) + \left(c_2 - a_2\right) - \left(c_3 - a_3\right) \right)^2  \\  \nonumber
4 (c_1 - a_1) (c_3 - a_3) & =  \;\left(\;\left(c_1 - a_1\right) - \left(c_2 - a_2\right) + \left(c_3 - a_3\right) \right)^2  \\
4 (c_2 - a_2) (c_3 - a_3) & =  \left(\shortminus\left(c_1 - a_1\right) + \left(c_2 - a_2\right) + \left(c_3 - a_3\right) \right)^2
\end{align}
defining  rational Pythagorean quadruples which can be normalized to define rational points on the unit sphere.
\newline
The Euclidean distance relationship becomes clear dividing both sides of $(16)$ by $2$ and taking the square root.
\newline
For the example $(m,n)=(2,1)$ we get the vectors 
$\mathbf{a}=\left(\frac{15}{2}, \frac{40}{3}, \frac{24}{5}\right)$ and  $\mathbf{c}=\left(\frac{353}{30}, \frac{881}{60}, \frac{337}{60}\right)$ \newline
for which the identity evaluates to $\left(\frac{193}{30}\right)^2= \left(\frac{24}{5}\right)^2+\left(\frac{56}{15}\right)^2+\left(\frac{21}{10}\right)^2$ 
and the 2-Norm gives us $\|c-a\|_2=\frac{1}{\sqrt{2}}\frac{193}{30}$.
\pagebreak
\subsection{The Trinity system}
The Trinity system is a three-dimensional rational differential geometric system originating from the identity
\begin{align}
(a_1^2 +a_2^2+a_3^2) = 2 (b_1^2 +b_2^2+b_3^2)=\frac{2}{3}(c_1^2 +c_2^2+c_3^2)=\frac{4\left(C^4 - (A B)^2\right)^2}{(A  B C)^2}
\end{align} 
relating the sum of the squares of the sides of the rational triples. Viewing the sides as cartesian coordinates  
\begin{align}&(x,y,z)_1=(a_1,a_2,a_3) &  &(x,y,z)_2=(b_1,b_2,b_3)&   & (x,y,z)_3=(c_1,c_2,c_3)&\end{align}
we obtain after substituting the ratio $\frac{m}{n}$ by the variable $t$, three normalized equations $S_i: x_i^2+y_i^2+z_i^2=r_i^2$ 
\small
\begin{eqnarray}
S_1: &(x,y,z,r)_1 &=\left(\frac{\left(t^4-1\right)^2}{t^8+14 t^4+1},\frac{4 t^2\left(t^2+1\right)^2}{t^8+14
   t^4+1},\frac{4 t^2 \left(t^2-1\right)^2}{t^8+14 t^4+1},1\right) \\
S_{2}:& (x,y,z,r)_2 &= \left(\frac{4 t^2 \left(t^4+1\right)}{t^8+14 t^4+1},\frac{\left(t^2-1\right)^2
   \left(t^4+6 t^2+1\right)}{2 \left(t^8+14 t^4+1\right)},-\frac{\left(t^2+1\right)^2
   \left(t^4-6 t^2+1\right)}{2 \left(t^8+14 t^4+1\right)},\sqrt{\frac{1}{2}}\right)    \\
S_{3}:& (x,y,z,r)_3 & =  \left(\frac{t^8+6 t^4+1}{t^8+14 t^4+1},\frac{t^8+4 t^6+22 t^4+4 t^2+1}{2 \left(t^8+14
   t^4+1\right)},\frac{t^8-4 t^6+22 t^4-4 t^2+1}{2 \left(t^8+14 t^4+1\right)},\sqrt{\frac{3}{2}}\right)
\end{eqnarray}
with the relationships 
\begin{align}
 x_1 + y_1 - z_1 &= 1  \;\;\; \;\;\;\;\;\;\;\;  &x_2 - y_2 - z_2 &= 0  \;\;\; \;\;\;\;\;\;\;\; &x_3 + y_3 + z_3 &= 2  \\
 x_1^2 + y_1^2 + z_1^2 &= 1\;\;\; \;\;\;\;\;\;\;\; &x_2^2 + y_2^2 + z_2^2 &= 1/2  \;\;\; \;\;\;\;\;\;\;\; & x_3^2 + y_3^2 + z_3^2 &= 3/2  \\
 x_1^2 + x_2^2 &= x_3^2\;\;\; \;\;\;\;\;\;\;\; &y_1^2 + y_2^2 &= y_3^2 \;\;\; \;\;\;\;\;\;\;\;& z_1^2 + z_2^2 &= z_3^2 
\end{align}
\begin{align}
 \frac{\partial ^n x_1}{\partial t^n} +  \frac{\partial ^n y_1}{\partial t^n} - \frac{\partial ^n z_1}{\partial t^n} &=0 
 &\frac{\partial ^n x_2}{\partial t^n} -  \frac{\partial ^n y_2}{\partial t^n} - \frac{\partial ^n z_2}{\partial t^n} &=0 
 &\frac{\partial ^n x_3}{\partial t^n} +  \frac{\partial ^n y_3}{\partial t^n} + \frac{\partial ^n z_3}{\partial t^n} &=0  & \forall \;n>0
\end{align}
\subsubsection{The circles }
The locus of the points $S_i:(x,y,z)_i$ on the spheres with radius $r_i=\{1,\sqrt{\frac{1}{2}},\sqrt{\frac{3}{2}}\}$ , \newline including their permutations over the 3 coordinates, define the circles $C_i(t)_{j}$ :
\newline
$C_1(t)_{1..8}$ in the planes $ \pm x \pm y\pm z =1$ with radius $r_{1,j}=\sqrt{\frac{2}{3}}$ and center at $(x,y,z)_{1,j}=(\pm\frac{1}{3},\pm\frac{1}{3},\pm\frac{1}{3})$.
\newline
$C_2(t)_{1..4}$ in the planes $\pm x \pm y + z=0$ with radius $r_{2,j}=\sqrt{\frac{1}{2}}$ and center at $(x,y,z)_{2,j}=(0,0,0)$.
\newline
$C_3(t)_{1..8}$ in the planes $\pm x \pm y \pm z =2$ with radius $r_{3,j}=\sqrt{\frac{1}{6}}$ and center at $(x,y,z)_{3,j}=(\pm\frac{2}{3},\pm\frac{2}{3},\pm\frac{2}{3})$.
\newline
The circles $C_1(t)_{1..8}$ also define the intersectionpaths between the sphere with radius $1$ centered at the origin and the $8$ spheres with radius $\sqrt{2}$ centered at $(\pm1,\pm1,\pm1)$ and
\newline
the cirlces $C_3(t)_{1..8}$ also define the intersectionpaths between the sphere with radius $\sqrt{\frac{3}{2}}$ centered at the origin and the $8$ spheres with radius $\frac{\sqrt{3}}{4}$ centered at $(\pm\frac{3}{4},\pm\frac{3}{4},\pm\frac{3}{4})$.
\newline
The figure below shows the circles $C_{1..3}$ and their basic trigonometric cartesian component functions.
\newline
\begin{figure}[hbt!]
\begin{tikzpicture}[scale=1.1,baseline={(current bounding box.center)}]
   \begin{axis}[
    		   view={-45}{-5},
                axis lines=center,axis on top,
                xlabel=$x$,ylabel=$y$,zlabel=$z$,
                xtick={5},ytick={5},ztick={5},
                no marks,axis equal,
                xmin=-1.2,xmax=1.2,ymin=-1.2,ymax=1.2,zmin=-1.2,zmax=1.2,
                enlargelimits={upper=0.1}]
	\addplot3+[blue,no markers,samples=100, samples y=0,domain=-pi:pi,variable=\t]
       ({(1/3-1/sqrt(3)*cos(\t r) -1/3*sin(\t r))},{(1/3+1/sqrt(3)*cos(\t r) -1/3*sin(\t r))},{(1/3+2/3*sin(\t r))});                                     
	\addplot3+[blue,no markers,samples=100, samples y=0,domain=-pi:pi,variable=\t]
       ({-(1/3-1/sqrt(3)*cos(\t r) -1/3*sin(\t r))},{(1/3+1/sqrt(3)*cos(\t r) -1/3*sin(\t r))},{(1/3+2/3*sin(\t r))});                                     
	\addplot3+[blue,no markers,samples=100, samples y=0,domain=-pi:pi,variable=\t]
       ({(1/3-1/sqrt(3)*cos(\t r) -1/3*sin(\t r))},{-(1/3+1/sqrt(3)*cos(\t r) -1/3*sin(\t r))},{(1/3+2/3*sin(\t r))});                                     
	\addplot3+[blue,no markers,samples=100, samples y=0,domain=-pi:pi,variable=\t]
       ({(1/3-1/sqrt(3)*cos(\t r) -1/3*sin(\t r))},{(1/3+1/sqrt(3)*cos(\t r) -1/3*sin(\t r))},{-(1/3+2/3*sin(\t r))});                                     
	\addplot3+[blue,no markers,samples=100, samples y=0,domain=-pi:pi,variable=\t]
       ({-(1/3-1/sqrt(3)*cos(\t r) -1/3*sin(\t r))},{-(1/3+1/sqrt(3)*cos(\t r) -1/3*sin(\t r))},{(1/3+2/3*sin(\t r))});                                     
	\addplot3+[blue,no markers,samples=100, samples y=0,domain=-pi:pi,variable=\t]
       ({-(1/3-1/sqrt(3)*cos(\t r) -1/3*sin(\t r))},{(1/3+1/sqrt(3)*cos(\t r) -1/3*sin(\t r))},{-(1/3+2/3*sin(\t r))});                                     
	\addplot3+[blue,no markers,samples=100, samples y=0,domain=-pi:pi,variable=\t]
       ({(1/3-1/sqrt(3)*cos(\t r) -1/3*sin(\t r))},{-(1/3+1/sqrt(3)*cos(\t r) -1/3*sin(\t r))},{-(1/3+2/3*sin(\t r))});                                     
	\addplot3+[blue,no markers,samples=100, samples y=0,domain=-pi:pi,variable=\t]
       ({-(1/3-1/sqrt(3)*cos(\t r) -1/3*sin(\t r))},{-(1/3+1/sqrt(3)*cos(\t r) -1/3*sin(\t r))},{-(1/3+2/3*sin(\t r))});                                     

	\addplot3+[green,no markers,samples=100, samples y=0,domain=-pi:pi,variable=\t]
	({(-1/2*cos(\t r)-1/(2*sqrt(3))*sin(\t r))},{(-1/2*cos(\t r)+1/(2*sqrt(3))*sin(\t r))},{(-1/sqrt(3)*sin(\t r))});
	\addplot3+[green,no markers,samples=100, samples y=0,domain=-pi:pi,variable=\t]
	({-(-1/2*cos(\t r)-1/(2*sqrt(3))*sin(\t r))},{(-1/2*cos(\t r)+1/(2*sqrt(3))*sin(\t r))},{(-1/sqrt(3)*sin(\t r))});
	\addplot3+[green,no markers,samples=100, samples y=0,domain=-pi:pi,variable=\t]
	({(-1/2*cos(\t r)-1/(2*sqrt(3))*sin(\t r))},{-(-1/2*cos(\t r)+1/(2*sqrt(3))*sin(\t r))},{(-1/sqrt(3)*sin(\t r))});
	\addplot3+[green,no markers,samples=100, samples y=0,domain=-pi:pi,variable=\t]
	({(-1/2*cos(\t r)-1/(2*sqrt(3))*sin(\t r))},{(-1/2*cos(\t r)+1/(2*sqrt(3))*sin(\t r))},{-(-1/sqrt(3)*sin(\t r))});

	\addplot3+[red,no markers,samples=100, samples y=0,domain=-pi:pi,variable=\t]
	({2/3-1/(2*sqrt(3))*cos(\t r) -1/6*sin(\t r)},{2/3+1/(2*sqrt(3))*cos(\t r) -1/6*sin(\t r)},{2/3+1/3*sin(\t r)});
	\addplot3+[red,no markers,samples=100, samples y=0,domain=-pi:pi,variable=\t]
	({-(2/3-1/(2*sqrt(3))*cos(\t r) -1/6*sin(\t r))},{(2/3+1/(2*sqrt(3))*cos(\t r) -1/6*sin(\t r))},{(2/3+1/3*sin(\t r))});
	\addplot3+[red,no markers,samples=100, samples y=0,domain=-pi:pi,variable=\t]
	({(2/3-1/(2*sqrt(3))*cos(\t r) -1/6*sin(\t r))},{-(2/3+1/(2*sqrt(3))*cos(\t r) -1/6*sin(\t r))},{(2/3+1/3*sin(\t r))});
	\addplot3+[red,no markers,samples=100, samples y=0,domain=-pi:pi,variable=\t]
	({(2/3-1/(2*sqrt(3))*cos(\t r) -1/6*sin(\t r))},{(2/3+1/(2*sqrt(3))*cos(\t r) -1/6*sin(\t r))},{-(2/3+1/3*sin(\t r))});

	\addplot3+[red,no markers,samples=100, samples y=0,domain=-pi:pi,variable=\t]
	({-(2/3-1/(2*sqrt(3))*cos(\t r) -1/6*sin(\t r))},{-(2/3+1/(2*sqrt(3))*cos(\t r) -1/6*sin(\t r))},{(2/3+1/3*sin(\t r))});
	\addplot3+[red,no markers,samples=100, samples y=0,domain=-pi:pi,variable=\t]
	({-(2/3-1/(2*sqrt(3))*cos(\t r) -1/6*sin(\t r))},{(2/3+1/(2*sqrt(3))*cos(\t r) -1/6*sin(\t r))},{-(2/3+1/3*sin(\t r))});
	\addplot3+[red,no markers,samples=100, samples y=0,domain=-pi:pi,variable=\t]
	({(2/3-1/(2*sqrt(3))*cos(\t r) -1/6*sin(\t r))},{-(2/3+1/(2*sqrt(3))*cos(\t r) -1/6*sin(\t r))},{-(2/3+1/3*sin(\t r))});
	\addplot3+[red,no markers,samples=100, samples y=0,domain=-pi:pi,variable=\t]
	({-(2/3-1/(2*sqrt(3))*cos(\t r) -1/6*sin(\t r))},{-(2/3+1/(2*sqrt(3))*cos(\t r) -1/6*sin(\t r))},{-(2/3+1/3*sin(\t r))});
                                      
  \end{axis}
\end{tikzpicture}
\begin{minipage}{0.6\textwidth}
\begin{align} \nonumber
\color{blue} C_1(t)= &\left(\frac{1}{3} -\frac{1}{\sqrt{3}}\cos(t)-\frac{1}{3}\sin(t)  ,\frac{1}{3} +\frac{1}{\sqrt{3}}\cos(t) -\frac{1}{3}\sin(t)  ,\frac{1}{3} +\frac{2}{3}\sin(t)\right) \\ \nonumber
\color{green} C_2(t)=  &\left(-\frac{1}{2}\cos(t)-\frac{1}{2\sqrt{3}}\sin(t)  ,-\frac{1}{2}\cos(t)+\frac{1}{2\sqrt{3}}\sin(t)  ,-\frac{1}{\sqrt{3}}\sin(t)\right) \\
\color{red} C_3(t)= &\left(\frac{2}{3} -\frac{1}{2\sqrt{3}} \cos(t)-\frac{1}{6} \sin(t)  ,\frac{2}{3} +\frac{1}{2\sqrt{3}} \cos(t) -\frac{1}{6}\sin(t)  ,\frac{2}{3} -\frac{1}{6}\sin(t)\right)
\end{align}
\end{minipage}\hspace{10cm}
\caption{The Trinity circles} 
\label{fig:F1}
\end{figure}
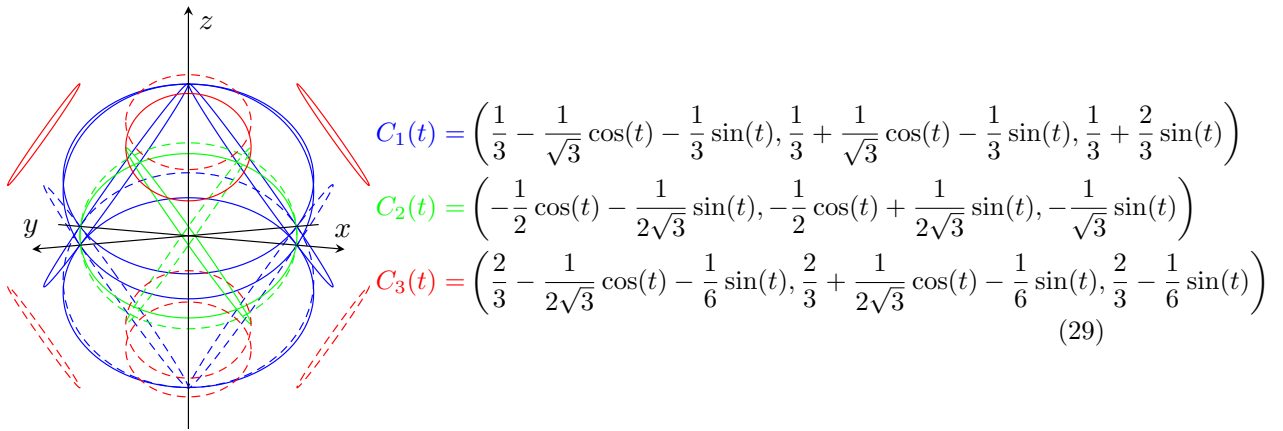
\pagebreak
\subsubsection{The vectors }
Defining the vectors $\mathbf {a},\mathbf {b},\mathbf {c}$ using the equations $(22,23,24)$ as follows 
\begin{align}
 \mathbf {a} &= (x,y,z)_1\;,&\mathbf {b} & = (x,-y,z)_2 \;, & \mathbf {c} & = (x,y,-z)_3& 
\end{align}
the Euclidean distances to the origin for the three vectors are 
\begin{align}
 &\left\|\mathbf {a} \right\|^2=1\;, &\;\;\; &\left\|\mathbf {b} \right\|^2=\frac{1}{2}\;, \;\;\; &\left\|\mathbf {c} \right\|^2&=\frac{3}{2}& 
\end{align}
Evaluating the dot products of the three vector pairs 
\begin{align}
& \mathbf {a}\cdot \mathbf {b}=0 \;, &\;\;\; &  \mathbf {b}\cdot \mathbf {c}=0 \;, &\;\;\; &\mathbf {a}\cdot \mathbf {c}=1 &
\end{align}
we find that the vector pairs $\mathbf {a},\mathbf {b}$ and $\mathbf {b},\mathbf {c}$ are orthogonal and between the vectors $\mathbf {a},\mathbf {c}$ we have a specific angle defined by 
\begin{equation}
\cos{\phi}=\frac{\mathbf {a} \cdot \mathbf{c}}{\left\|\mathbf {a} \right\|\left\|\mathbf {c} \right\|}= \sqrt{\frac{2}{3}}
\;\;\;,\;\;\;\phi = \arccos{\sqrt{\frac{2}{3}}}\end{equation}
which is an angle of $35.2643...$°.
\newline
An interesting fact is that the angle bewteen the cross product of the vector pair $\mathbf {a},\mathbf {b}$ and the vector $\mathbf {c}$ as well as for the  cross product of the vector pair $\mathbf {b},\mathbf {c}$ and the vector $\mathbf {a}$  are equal and can be evaluated by 
\begin{equation}
\cos{\theta}=\frac{(\mathbf {a} \times \mathbf {b}) \cdot \mathbf{c}}{\left\|\mathbf {a} \times \mathbf {b} \right\|\left\|\mathbf {c} \right\|}= \frac{(\mathbf {b} \times \mathbf {c}) \cdot \mathbf{a}}{\left\|\mathbf {b} \times \mathbf {c} \right\|\left\|\mathbf {a} \right\|}=\sqrt{\frac{1}{3}} \;\;\;,\;\;\;\theta = \arccos{\sqrt{\frac{1}{3}}}
\end{equation}
which is an angle of $54.7356...$°.
\newline
\linebreak
Taking the dot products of the $n$-th derivative of the vectors denoted by $\frac{\partial ^n }{\partial t^n} $ we obtain  
\begin{align}
  &\frac{\partial ^n \mathbf {a}}{\partial t^n} \cdot \frac{\partial ^n \mathbf {b}}{\partial t^n}  =0 \;, &\;\;\; & \frac{\partial ^n\mathbf {b} }{\partial t^n} \cdot \frac{\partial ^n\mathbf {c} }{\partial t^n} =0 \;, &\;\;\; & \frac{\partial ^n \mathbf {a}}{\partial t^n}  \cdot \frac{\partial ^n \mathbf {c}}{\partial t^n} = \frac{1}{2} \left\|\frac{\partial ^n \mathbf {a}}{\partial t^n}  \right\|^2 = \frac{2}{3} \left\|\frac{\partial ^n \mathbf {b}}{\partial t^n}  \right\|^2 = 2 \left\|\frac{\partial ^n \mathbf {c} }{\partial t^n}  \right\|^2&
\end{align}
resulting in the norm equations for the derived spheres and their rational points
\begin{align} &
\left\|\frac{\partial ^n \mathbf {a} }{\partial t^n}  \right\|^2=2 \frac{\partial ^n  \mathbf {a} }{\partial t^n}\cdot \frac{\partial ^n  \mathbf {c} }{\partial t^n}\;, \;\;\; & &
\left\|\frac{\partial ^n \mathbf {b} }{\partial t^n}\right\|^2=\frac{3}{2}\frac{\partial ^n  \mathbf {a} }{\partial t^n}\cdot \frac{\partial ^n  \mathbf {c} }{\partial t^n} \;, \;\;\; &
\left\|\frac{\partial ^n  \mathbf {c}}{\partial t^n}  \right\|^2=\frac{1}{2} \frac{\partial ^n  \mathbf {a} }{\partial t^n}\cdot \frac{\partial ^n  \mathbf {c} }{\partial t^n} &
 \end{align}
 \newline
The scalar triple product of $\mathbf {a},\mathbf {b},\mathbf {c}$ is equal to $\frac{1}{2}$ and remains unchanged under a circular shift of its operands
\begin{gather}
 \mathbf {a}\cdot\left(\mathbf {b}\times \mathbf {c}\right)=\; \mathbf {b}\cdot\left(\mathbf {c}  \times \mathbf {a}\right)=\; \mathbf {c} \cdot\left( \mathbf {a}\times\ \mathbf {b}\right)=\;\;\frac{1}{2}  
 \end{gather}
It also represents geometrically the volume of the paralellepipeds created by the vector vertices. One sixth of its value, which is $\frac{1}{12}$, represents the volume of the tetrahedrons defined by the three vector vertices and the origin.
\newline
\linebreak
The vector triple product defined as the cross product of $\mathbf {a}$ or $\mathbf {c}$  with the cross product of $\mathbf {b}$ and $\mathbf {c}$ or $\mathbf {a}$ can be simplified due to $(32)$ and the relationship $\mathbf {a}\times \left( \mathbf {b} \times \mathbf {c}\right) = \left(\mathbf {a}\cdot \mathbf {c}\right)\mathbf {b} -  \left(\mathbf {a}\cdot \mathbf {b}\right)\mathbf {c} $
to 
\begin{gather}
 \mathbf {a} \times\left(\mathbf {b}  \times \mathbf {c}\right)=\mathbf {c} \times\left(\mathbf {b}\times \mathbf {a}\right)  = \mathbf {c} \times \mathbf {a }= \mathbf {b} 
\end{gather}
The vector triple product of $\mathbf {b}$ with $\mathbf {a} \times \mathbf {c}$ is the zero vector because the $\mathbf {a}$ and $\mathbf {c}$ are orthogonal to $\mathbf {b}$
\begin{gather}
 \mathbf {b} \times\left(\mathbf {a} \times \mathbf {c}\right)= \mathbf{0}  
\end{gather}
Other symmetry relationships involving cross and dot products of the vector derivatives are :
\begin{gather}
  \frac{\partial ^n \mathbf {a}}{\partial t^n} \times \frac{\partial ^n \mathbf {c}}{\partial t^n} =\mathbf{0} \;\; \forall \;\; n>0 \\
  2\frac{\partial ^n\mathbf {b} }{\partial t^n} \cdot \frac{\partial ^m\mathbf {c} }{\partial t^m} = \frac{\partial ^n \mathbf {b}}{\partial t^n}  \cdot \frac{\partial ^m \mathbf {a}}{\partial t^m} 
\\
  2\frac{\partial ^n\mathbf {b} }{\partial t^n} \times \frac{\partial ^m\mathbf {c} }{\partial t^m} = \frac{\partial ^n \mathbf {b}}{\partial t^n}  \times \frac{\partial ^m \mathbf {a}}{\partial t^m} 
 \\ 
 3  \frac{\partial ^n \mathbf {a}}{\partial t^n} \cdot   \frac{\partial ^m \mathbf {a}}{\partial t^m} =  4  \frac{\partial ^n \mathbf {b}}{\partial t^n} \cdot   \frac{\partial ^m \mathbf {b}}{\partial t^m} =  12  \frac{\partial ^n \mathbf {c}}{\partial t^n} \cdot   \frac{\partial ^m \mathbf {c}}{\partial t^m}
\\ 3  \frac{\partial ^n \mathbf {a}}{\partial t^n} \times   \frac{\partial ^m \mathbf {a}}{\partial t^m} =  4  \frac{\partial ^n \mathbf {b}}{\partial t^n} \times   \frac{\partial ^m \mathbf {b}}{\partial t^m} =  12  \frac{\partial ^n \mathbf {c}}{\partial t^n} \times   \frac{\partial ^m \mathbf {c}}{\partial t^m} 
\end{gather}
\begin{gather}
 \frac{\partial ^n \mathbf {a}}{\partial t^n} \cdot \frac{\partial ^m \mathbf {c}}{\partial t^m}(-1,-1,1) 
= 2  \frac{\partial ^n \mathbf {b}}{\partial t^n} \times \frac{\partial ^m \mathbf {c}}{\partial t^m} 
=  \frac{\partial ^n \mathbf {b}}{\partial t^n} \times \frac{\partial ^m \mathbf {a}}{\partial t^m} 
\end{gather}
\begin{gather}
3 \frac{\partial ^n \mathbf {a}}{\partial t^n} \times \frac{\partial ^m \mathbf {c}}{\partial t^m} 
= 2 \frac{\partial ^n \mathbf {b}}{\partial t^n} \cdot \frac{\partial ^m \mathbf {c}}{\partial t^m} (1,1,-1) 
= \frac{\partial ^n \mathbf {b}}{\partial t^n} \cdot \frac{\partial ^m \mathbf {a}}{\partial t^m} (1,1,-1)\;\; \forall \;\;\{m,n\}> 0 
\end{gather}
The equations and examples can be verified in the Sagemath notebook $\texttt{02\_Beyond\_Pythagoras.ipynb}$.
\pagebreak
\section{The conjugate conics}
\label{section:Ellipse}
We define a geometric method relating the congruent number problem to two conjugate conics, an ellipse $e$ and a degenerate hyperbola $h$ touching eachother. By this method congruent numbers and polynomials can be defined from intersections of a line and an ellipse. 
To see this relationship a right triangle with sides $a,b,c\in \mathbb{Q}$ for a congruent number $N$ is defined using the conjugate conics $e,h$ and the variables $f_1,f_2$ by
\small
\begin{equation}(e,h)=\left(\sqrt{N f_1^2 f_2^2- \frac{1}{4}(N f_1^2-f_2^2)^2},\sqrt{N f_1^2 f_2^2+ \frac{1}{4}(N f_1^2-f_2^2)^2}\right)\end{equation}
$$(a,b,c) = \left(\frac{2\; e\; h}{f_1 f_2 (N f_1^2-f_2^2)}, \frac{N f_1 f_2 (N f_1^2-f_2^2) }{ e\; h },\sqrt{\left(\frac{2\; e\; h}{f_1 f_2 (N f_1^2-f_2^2)}\right)^2+\left(\frac{ N f_1 f_2 (N f_1^2-f_2^2)}{e\; h}\right)^2 } \right)$$
for which $h$ reduces to $\frac{1}{2}(N f_1^2+f_2^2)$ and $(a,b,c)$ simplifies to
\begin{equation}(a,b,c) =
\left(\frac{ e\; (N f_1^2+f_2^2)}{f_1 f_2 (N f_1^2-f_2^2)}, \frac{2 N f_1 f_2 (N f_1^2-f_2^2) }{e\; (N f_1^2+f_2^2) },\frac{
4 N f_1^2 f_2^2\left(3 N f_1^2 f_2^2-(N f_1^2-f_2^2)^2 \right)+(N^2 f_1^4+f_2^4)^2}{4\; e\; f_1 f_2 (N^2 f_1^4 -f_2^4) }\right)
\end{equation}
\newline
For certain $N,f_1,f_2$ , $e$ becomes rational and also the triple $(a,b,c)$. Looking at $N$ as the independent variable $x$, the equation $e(x)$ is an ellipse and all its rational points can be found by intersecting $e(x)$ with a line $y(x)=t(x-x_0)+y_0$ going through a rational point $Q:(x_0,y_0)$ on the ellipse.
\newline
\linebreak
The figure below shows the translated ellipse $e(x)$ touching the degenerate hyperbola $h(x)$ at the point $Q:(x_0,y_0)$.
\newline
The line $y(x)$ intersects the ellipse in $Q$ and a second rational point $P:(x_t,e_t)$. 
\newline 
From $x_t$ we obtain a primitive congruent number $N$ using a procedure to reduce and or raise $x_t\in \mathbb{Q}$ to $N\in \mathbb{Z}$.
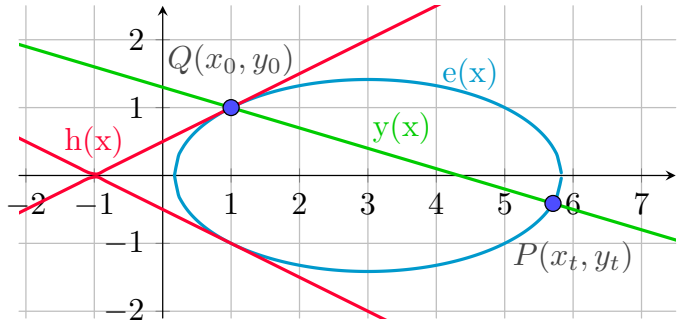
\begin{figure}[hbt!]
\centering
\begin{tikzpicture}[scale=1.2, line cap=round,line join=round,>=triangle 45,x=1cm,y=1cm]
\begin{axis}[
x=0.75cm,y=0.75cm,
axis lines=middle,
ymajorgrids=true,
xmajorgrids=true,
xmin=-2.1,
xmax=7.5,
ymin=-2.1,
ymax=2.5,
xtick={-2,-1,...,7},
ytick={-2,-1,...,2},]
\clip(-2.6541760682099147,-6.42426466103129) rectangle (10.277881323330046,5.980780987279769);

\draw[line width=1pt,color=qqzzcc,smooth,samples=100,domain=0.17157287586052444:5.828427030452535] plot(\x,{sqrt((\x)-1/4*(-1+(\x))^(2))});
\draw[line width=1pt,color=ffqqtt,smooth,samples=100,domain=-2.6541760682099147:10.277881323330046] plot(\x,{sqrt((\x)+1/4*(-1+(\x))^(2))});
\draw[line width=1pt,color=qqccqq,smooth,samples=100,domain=-2.6541760682099147:10.277881323330046] plot(\x,{1-0.3*(-1+(\x))});
\draw[line width=1pt,color=qqzzcc,smooth,samples=100,domain=0.17157287586052444:5.828427030452535] plot(\x,{0-sqrt((\x)-1/4*(-1+(\x))^(2))});
\draw[line width=1pt,color=ffqqtt,smooth,samples=100,domain=-2.6541760682099147:10.277881323330046] plot(\x,{0-sqrt((\x)+1/4*(-1+(\x))^(2))});
\begin{scriptsize}
\draw[color=qqzzcc] (4.5,1.5) node {e(x)};
\draw[color=ffqqtt] (-1.0,0.5) node {h(x)};
\draw[color=qqccqq] (3.5,0.7) node {y(x)};
\draw [fill=ududff] (1.0,1.0) circle (2.5pt);
\draw[color=uuuuuu] (1.0,1.7) node {$Q(x_0,y_0)$};
\draw [fill=ududff] (5.705882352940685,-0.4117647058822054) circle (2.5pt);
\draw[color=uuuuuu] (6.0,-1.2) node {$P(x_t,y_t)$};
\end{scriptsize}
\end{axis}
\end{tikzpicture}
\caption{Congruent numbers as intersections of a line and an ellipse} 
\label{fig:F2}
\end{figure}
The process of reducing due to square factors of $x_t$ is to divide the sides $a,b$ by the common factor reducing the area by the square of this factor.
\newline
The process of raising due to the denominator of $x_t$ is to multiply $a,b$ by this denominator. 
\newline
Both procedures are applied in the example given in the section  $\ref{section:EP}$.
\newline
\linebreak
The congruent number elliptic curve $ E_N: y^2= x^3-N^2 x$ has the following 2 points of infinite order
\begin{eqnarray}
P_1:(x,y)_1 &= &\left(-\frac{\left(N f_1^2 -f_2^2\right)^2}{4 f_1^2 f_2^2},\frac{\left(N f_1^2-f_2^2\right)^5-16 N^2 f_1^4 f_2^4 \left(N f_1^2 -f_2^2\right)}{32 e h f_1^3  f_2^3 }\right) \\
P_2:(x,y)_2 &=& \left(\frac{4N^2 f_1^2 f_2^2}{(N f_1^2 -f_2^2)^2},\frac{16 N^4 f_1^5 f_2^5-N^2 f_1 f_2 \left(N f_1^2-f_2^2\right)^4}{2 e h \left(N f_1^2-f_2^2\right)^3}\right) 
\end{eqnarray}
The equations are also valid for certain congruent numbers which are prime or an almost prime by adjoining $f_2$.
\newline
From the footprint equations defined in section $\ref{section:FP}$ we can verify the following congruence conditions for primes $p$ and the term by which $f_2$ needs to be adjoined  
$$ \begin{array}{ccccccc}
 N= &p & \forall & p \equiv\; 5 \mod 8 &   f_1 \in \mathbb{Z} ,& f_2 \in  \mathbb{Z} \;\;\;\;\;\;\;\;\;\;\; &\\
 N= &p & \forall & p \equiv\; 7 \mod 8 &   f_1 \in \mathbb{Z} , & f_2 \in \mathbb{Z}[\sqrt{N}]\;\;\;& \\
 N= &2p & \forall & p \equiv\; 7 \mod 8 &  f_1 \in \mathbb{Z} , & f_2 \in \mathbb{Z}[\sqrt{2N}]& \\
\end{array}$$ 
In this case $f_1,f_2$ are equal to the footprint solutions $m,n$ for certain primes $p$ listed in the tables $\ref{table:TI}\;,\ref{table:TII}\;,\ref{table:TIII}$.
\newline
\subsection{A famous example}
The famous example due to Don Zagier who found the smallest rational triangle $(a,b,c)$ for the congruent number $N=157$ can be verified  evaluating $(48)$ for $(f_1,f_2)=(87005 ,610961)$ giving the triple 
\begin{align*}
(a,b,c) =\Bigg( &   \frac{411340519227716149383203}{21666555693714761309610} ,  \frac{6803298487826435051217540}{411340519227716149383203}& \\
 ,&\;\; \frac{224403517704336969924557513090674863160948472041}{8912332268928859588025535178967163570016480830}\Bigg)
\end{align*}   
and evaluating $(49,50)$ we obtain the rational points for the elliptic curve $E: y^2 = x^3-157^2x$ 
\begin{align*}
P_1&:\left(-\frac{166136231668185267540804}{2825630694251145858025},-\frac{1
   67661624456834335404812111469782006}{150201095200135518108761470235125}\right) &\\
P_2&:\left(\frac{69648970982596494254458225}{166136231668185267540804},\frac
   {538962435089604615078004307258785218335}{67716816556077455999228495435742408}\right)&
\end{align*}   
\subsection{The intersection example}
\label{section:EP}
By the line and  ellipse intersection method we prove that $$N(t)=\left(4 t^2+1\right) \left(4 t^2-8 t+5\right)$$  
is a congruent number polynomial.
\newline
\linebreak
First we take $(47)$ and change the variables $N=x$ and $(f_1,f_2)=(1,f)$ , such that the ellipse and hyperbola $$(e,h)=\left(\sqrt{f^2 x-\frac{1}{4} \left(x-f^2\right)^2},\frac{1}{2} \left(f^2+x\right)\right)$$
are conjugate conics touching at the point $Q=(f^2,f^2)$. To find a rational right triangle $(a,b,c)$ as defined by $(48)$ we only have to find the rational points on the ellipse because the hyperbola is degenerate and reduces to a set of lines.
Intersecting the ellipse $e(x)$ with a line $y(x)=t(x-f^2)+f^2$ going through $Q$ gives the parameterization of the rational points
$$P:(x,e)=\left(\frac{f^2 \left(4 t^2-8 t+5\right)}{4 t^2+1},\frac{f^2 \left(-4 t^2+4 t+1\right)}{4 t^2+1}\right)$$
We notice that we can reduce the sides $a,b$ by $f$ due to the square factor and raise the sides by $4 t^2+1$ due to the denominator.
Applying first the intersection solution $(x,e)$ to the rational right triangle $(a,b,c)$ defined by $(48)$ and secondly raising the sides by $4 t^2+1$ we obtain the right triangle definition   
\begin{align*}
(a,b,c)=&\Bigg( \frac{-16 t^4+32 t^3-24 t^2+8 t+3}{2-4 t} ,\frac{4 \left(32 t^5-80 t^4+80 t^3-40 t^2+18 t-5\right)}{16 t^4-32 t^3+24 t^2-8 t-3} & \\ 
,&\frac{256 t^8-1024 t^7+1792 t^6-1792 t^5+1504 t^4-1216 t^3+688 t^2-208 t+41}{64 t^5-160 t^4+160 t^3-80 t^2+4 t+6}\Bigg)&
\end{align*}
completing the proof for the congruent number polynomial    
$$N(t)=\frac{a b}{2}=\frac{-4 \left(32 t^5-80 t^4+80 t^3-40 t^2+18 t-5\right)}{2(2-4 t)} =\left(4 t^2+1\right) \left(4 t^2-8 t+5\right)$$
The points on the elliptic curve (49,50) are then :
\begin{align*}
P_1: & \left(-4 (2 t -1)^2, 2 (2 t-1) (4 t^2-4 t -1) (4 t^2-4 t+3) \right) \\
P_2: & \left( \frac{(4 t^2+1)^2 ( 4 t^2-8t+5)^2}{ 4 (2 t-1)^2}, \frac{(4 t^2+1)^2 (4 t^2-8t+5)^2 (4 t^2-4t-1) (4 t^2-4t+3)}{8 (2 t-1)^3}\right) 
\end{align*}
As an example for $N(3)=629$ we obtain the rational triangle  
 $$(a,b,c)=\left(\frac{621}{10},\frac{12580}{621},\frac{405641}{6210}\right)$$
 and the two points on the elliptic curve $ E: y^2= x^3-629^2 x$
 $$ P_1:(-100,6210) \;\;\;,\;\;\; P_2:\left(\frac{395641}{100},\frac{245693061}{1000} \right)$$
 \pagebreak
\subsection{The lattice example}
A more general example follows from the lattice points of an ellipse.
\newline
Assigning $(N,f_1,f_2)=(x, 1,m^2+ n^2)\;\; \forall \;\;{m,n} \in \mathbb{Z}$  
we obtain for the ellipse and the degenerate hyperbola
\begin{eqnarray} 
(e,h)&=& \left(\sqrt{N f_1^2 f_2^2 -\frac{1}{4} \left(N f_1^2-f_2^2\right)^2},\sqrt{N f_1^2 f_2^2 +\frac{1}{4} \left(N f_1^2-f_2^2\right)^2}\right)\\
 & = & \left(\sqrt{x \left(m^2+n^2\right)^2-\frac{1}{4}\left(\left(m^2+n^2\right)^2-x\right)^2},\frac{1}{2}\left(\left(m^2+n^2\right)^2+x\right)\right) 
\end{eqnarray} 
the following lattice points for the ellipse 
\begin{eqnarray} \nonumber
& P_1:(x,e)_1 = &\left((m^2 + n^2) (m^2 + 4 m n + 5 n^2)  \;,\; (m^2 + n^2)(m^2 + 2mn - n^2)  \right) \\ \nonumber
& P_2:(x,e)_2 = &\left((m^2 + n^2) (m^2 - 4 m n +  5 n^2)   \;,\; (m^2 + n^2)(m^2 - 2mn - n^2)  \right) \\ \nonumber
& P_3:(x,e)_3 = &\left((m^2 + n^2) (5 m^2 - 4 m n + n^2)    \;,\; (m^2 + n^2)(m^2 + 2mn - n^2) \right) \\
& P_4:(x,e)_4 = &\left((m^2 + n^2) (5 m^2 +  4 m n + n^2)  \;,\; (m^2 + n^2)(m^2 - 2mn - n^2) \right)
\end{eqnarray} 
The $x_i$ coordinates represent congruent numbers, except when $x_i$ is a square, due to the right triangles  
\begin{align*}
(a,b,c)_1=& \Bigg(\frac{\left(m^2+2 m n-n^2\right) \left(m^2+2 m n+3 n^2\right)}{2 n (m+n)},\frac{4 n (m+n)
   \left(m^2+n^2\right) \left(m^2+4 m n+5 n^2\right)}{\left(m^2+2 m n-n^2\right) \left(m^2+2 m n+3 n^2\right)}\\
   ,&\frac{m^8+8 m^7 n+28 m^6 n^2+56 m^5 n^3+94 m^4 n^4+152 m^3 n^5+172 m^2 n^6+104 m n^7+41
   n^8}{2 n (m+n) \left(m^2+2 m n-n^2\right) \left(m^2+2 m n+3 n^2\right)}\Bigg)\\
(a,b,c)_2=& \Bigg(\frac{\left(m^2-2 m n-n^2\right) \left(m^2-2 m n+3 n^2\right)}{2 n (m-n)},\frac{4 n (m-n) \left(m^2+n^2\right) \left(m^2-4 m n+5 n^2\right)}{\left(m^2-2 m n-n^2\right)\left(m^2-2 m n+3 n^2\right)}\\
   ,&\frac{m^8-8 m^7 n+28 m^6 n^2-56 m^5 n^3+94 m^4 n^4-152 m^3 n^5+172 m^2 n^6-104 m n^7+41 n^8}{2 m^5 n-10 m^4
   n^2+20 m^3 n^3-20 m^2 n^4+2 m n^5+6 n^6}\Bigg)\\
(a,b,c)_3=& \Bigg(\frac{\left(m^2+2 m n-n^2\right) \left(3 m^2-2 m n+n^2\right)}{2 m (m-n)},\frac{4 m (m-n) \left(m^2+n^2\right) \left(5 m^2-4 m n+n^2\right)}{\left(m^2+2 m n-n^2\right) \left(3m^2-2 m n+n^2\right)} \\
   ,&\frac{41 m^8-104 m^7 n+172 m^6 n^2-152 m^5 n^3+94 m^4 n^4-56 m^3 n^5+28 m^2 n^6-8 m n^7+n^8}{6 m^6+2 m^5 n-20
   m^4 n^2+20 m^3 n^3-10 m^2 n^4+2 m n^5}\Bigg)\\
(a,b,c)_4=& \Bigg(-\frac{\left(m^2-2 m n-n^2\right) \left(3 m^2+2 m n+n^2\right)}{2 m (m+n)},-\frac{4 m (m+n)
   \left(m^2+n^2\right) \left(5 m^2+4 m n+n^2\right)}{\left(m^2-2 m n-n^2\right) \left(3 m^2+2 m n+n^2\right)} \\
   ,&\frac{41 m^8+104 m^7 n+172 m^6 n^2+152 m^5 n^3+94 m^4 n^4+56 m^3 n^5+28 m^2 n^6+8 m n^7+n^8}{2 m \left(-3
   m^5+m^4 n+10 m^3 n^2+10 m^2 n^3+5 m n^4+n^5\right)}\Bigg)
\end{align*}
for $m,n$ such that the denominator of $a_{1,2,3,4}\ne 0$.
\newline
\linebreak
Intersecting the ellipse by the lines $y_i = t (x - x_i) + e_i$  we obtain a second rational point $P_{i2}$. 
\newline
Raising the $X$-axis intersection equations by $4 t^2+1$ we obtain the congruent number equations :
\begin{eqnarray}\nonumber
&N_{12}&=\left(4 t^2+1\right)\left(m^2+n^2\right) \left(m^2 \left(4 t(t-2)+5\right)+4 m n (4 t (t-1) -1)+n^2 (4 t (5 t+2)+1)\right)\\ \nonumber
&N_{22}&=\left(4 t^2+1\right)\left(m^2+n^2\right) \left(m^2 \left(4 t (t+2)+5\right)-4 m n (4 t (t+1)-1)+n^2 (4 t (5 t-2)+1)\right)\\ \nonumber
&N_{32}&=\left(4 t^2+1\right)\left(m^2+n^2\right) \left(m^2 \left(4 t (5 t-2)+1\right)-4 m n (4 t (t+1)-1)+n^2 (4 t (t+2)+5)\right)\\ 
&N_{42}&=\left(4 t^2+1\right)\left(m^2+n^2\right) \left(m^2 \left(4 t (5 t+2)+1\right)+4 m n (4 (t-1) t-1)+n^2 (4 (t-2) t+5)\right)
\end{eqnarray}
Evaluating for $(m,n,t)=(1,2,3)$ we obtain the congruent numbers and the right triangles
\begin{align*}
 N_{12} = 188885 \;\;\;&,\;\;\;  (a,b,c)_{12} = \left( \frac{71757}{418} , \frac{157907860}{71757} , \frac{66206019401}{29994426 }\right)   \\
 N_{22} = 58645  \;\;\;\;\;&,\;\;\;(a,b,c)_{22} = \left( \frac{58483}{66} , \frac{7741140}{58483} , \frac{3458210761}{3859878}\right)  \\
N_{32} = 7585  \;\;\;\;\;\;&,\;\;\;(a,b,c)_{32} = \left( \frac{ -5537}{72} , \frac{-1092240}{5537} , \frac{-84406081}{398664}\right) \\
N_{42} = 84545  \;\;\;\;\; &,\;\;\;(a,b,c)_{42} = \left( \frac{ 82497}{136} , \frac{22996240}{82497} , \frac{7489959041}{11219592 }\right) 
\end{align*}
\newline
The equations and examples for this section can be found in the Sagemath notebook $\texttt{03\_Conjugate\_conics.ipynb}$.
\pagebreak
\section{The twin hyperbolas}
\label{section:TH}
Another example of two combined conic sections related to the congruent number problem follows from the term under the square root $\sqrt{\frac{N}{2}\left((m-n)^2 +2m^2\right)\left((m+n)^2+2n^2\right)}$ of the footprint equation $(78)$ of section $\ref{section:FP}$. Replacing $m$ by $n x$ we get $n\sqrt{\frac{N}{2}\left(1-2x+3x^2\right)\left(3+2x+x^2\right)}$ obtaining the definition for the twin hyperbolas.
\newline
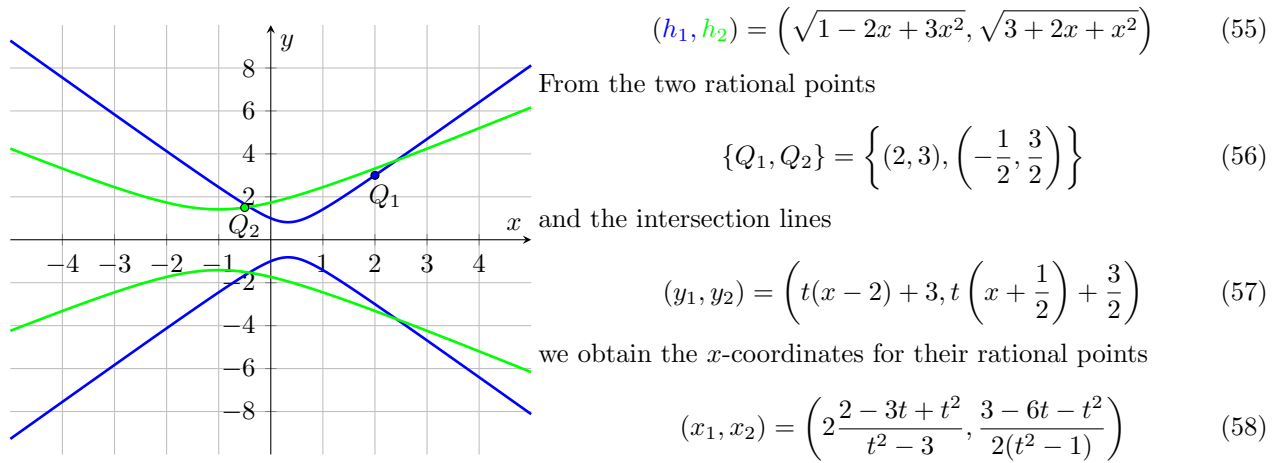
\begin{figure}[hbt!]
\begin{tikzpicture}[scale=1,baseline={(current bounding box.center)}]
\begin{axis}[
        xmin=-5,
        xmax=5,
        ymin=-10,
        ymax= 10,
        xlabel={$x$},
        ylabel={$y$},
        axis lines=middle,
        ymajorgrids=true,
        xmajorgrids=true,
        samples=100,
        smooth,
        xtick={-4,-3,...,4},
        ytick={-8,-6,...,8},
        clip=false,
    ]
\draw[line width=1pt,color=blue,smooth,samples=100,domain=-5:5] plot(\x, { sqrt(1-2*(\x)+3*(\x)^(2))});
\draw[line width=1pt,color=blue,smooth,samples=100,domain=-5:5] plot(\x, {-sqrt(1-2*(\x)+3*(\x)^(2))});
\draw[line width=1pt,color=green,smooth,samples=100,domain=-5:5] plot(\x, { sqrt(3+2*(\x)+(\x)^(2))});
\draw[line width=1pt,color=green,smooth,samples=100,domain=-5:5] plot(\x, {-sqrt(3+2*(\x)+(\x)^(2))});

\begin{scriptsize}
\draw[color= black]   (1.7,2) node[anchor = west] {$Q_1$};
\draw[color= black]   (-1,0.6) node[anchor = west] {$Q_2$};
\draw [fill= blue] (2,3) circle (1.5pt);
\draw [fill= green] (-0.5,1.5) circle (1.5pt);
\end{scriptsize}

\end{axis}
\end{tikzpicture}
\begin{minipage}{0.6\textwidth}
\begin{equation}
( \color{blue} h_1,\color{green}h_2 \color{black})=  \left(\sqrt{1-2x+3x^2},\sqrt{3+2x+x^2} \right)  
\end{equation}
From the two rational points 
\newline
\begin{equation}
\{Q_{1},Q_{2}\}=\left\{(2,3),\left(-\frac{1}{2},\frac{3}{2}\right)\right\}
\end{equation}
and the intersection lines 
\newline
\begin{equation}
(y_1,y_2)= \left( t (x-2) +3,  t \left(x+\frac{1}{2}\right) +\frac{3}{2}\right)
\end{equation}
we obtain the $x$-coordinates for their rational points 
\begin{equation}
(x_1,x_2) = \left( 2 \frac{2 - 3 t + t^2}{t^2-3} , \frac{3 - 6 t - t^2}{2(t^2-1)}\right)
\end{equation}
\end{minipage}\hspace{10cm}
\caption{The twin hyperbolas} 
\label{fig:F3}
\end{figure}
\newline
Evaluating the product $H(x)=2 h_1(x)^2 h_2(x)^2$ for the $x$-coordinates we get in function of $t$
\begin{equation} H(x_1)=
\left(\frac{9-10 t +3 t^2}{\left(t^2-3\right)^2}\right)^2  2  \left(11 t^4-36 t^3+30 t^2-12  t+19\right)
\end{equation}
\begin{equation} H(x_2)=
\left(\frac{3-2 t +3 t^2}{4\left(t^2-1\right)^2}\right)^2 \frac{\left(11 t^4+60 t^3+66 t^2-132  t+43\right)}{2}
\end{equation}
Then the square free part of $H(x) $ is related to a congruent number polynomial $N$ in function of $t$.
\newline
Taking the squarefree part from $(59)$ and four times the squarefree part of $(60)$ we obtain the two congruent number polynomials
\begin{equation} 
(N_1,N_2) = \left(  2 \left(11 t^4-36 t^3+30 t^2-12  t+19\right) , 2\left(11 t^4+60 t^3+66 t^2-132  t+43\right)\right)
\end{equation}
representing the areas for the rational right triangles with sides  
\begin{align*} (a,b)_1= & \Bigg(
-\frac{\left(3 t^2-10 t+9\right) \left(t^4+12 t^3-62 t^2+84 t-31\right) \left(11 t^4-36 t^3+30 t^2-12  t+19\right)}{\left(t^2-2 t-1\right) \left(7 t^2-22 t+17\right) \left(5 t^4-24 t^3+46 t^2-48 t+25\right)}
,\\ &
-\frac{4 \left(t^2-2 t-1\right) \left(7 t^2-22 t+17\right) \left(5 t^4-24 t^3+46 t^2-48 t+25\right)}{\left(3  t^2-10 t+9\right) \left(t^4+12 t^3-62 t^2+84 t-31\right)}
\Bigg)
\end{align*}
\begin{align*} (a,b)_2= & \Bigg(
-\frac{\left(3 t^2-2 t+3\right) \left(t^4+36 t^3+22 t^2-60 t+17\right) \left(11 t^4+60 t^3+66 t^2-132  t+43\right)}{\left(t^2+2 t-7\right) \left(7 t^2-2 t-1\right) \left(5 t^4+12 t^3+22 t^2-36 t+13\right)} 
,\\ &
-\frac{4 \left(t^2+2 t-7\right) \left(7 t^2-2 t-1\right) \left(5 t^4+12 t^3+22 t^2-36 t+13\right)}{\left(3   t^2-2 t+3\right) \left(t^4+36 t^3+22 t^2-60 t+17\right)}
\Bigg)
\end{align*}
For example for $t=10$ we obtain the congruent numbers and the right triangles
\begin{align*}
N_1=153798 \;\;\;,\;\;\; (a,b,c)_1 &= \left(-\frac{266938037619}{1183583135}, -\frac{4734332540}{ 3471281}, \frac{5679574272052285061}{4108549648445935}\right)  \\
N_2= 350646 \;\;\;,\;\;\; (a,b,c)_2 &= \left(-\frac{2362584547353}{4899249131}, -\frac{19596996524}{  13475611}, \frac{101151574309748379365}{66020375481444041}\right) 
\end{align*}
The equations and examples can be verified in the Sagemath notebook $\texttt{04\_Twin\_hyperbolas.ipynb}$.
\pagebreak
\section{The Cassini ovals}
The general equation for the Cassini oval,  with foci at $(a', 0)$ and $(-a', 0)$, is a quartic curve
$$ \left((x-a')^{2}+y^{2}\right)\left((x+a')^{2}+y^{2}\right)=b'^{4}$$
which can be written as a difference of two squares 
$ (y^2 + x^2+ a'^2)^2 - (2 a' x)^2 = b'^4 $ or as the function
\begin{equation}
y(x)=\sqrt{\sqrt{4 a'^2 x^2+b'^4}-x^2-a'^2}
\end{equation}
For $a'=1$ the Cassini oval intersects the coordinate axis depending on $b'$ at the points
$$\{(\shortminus c_4,0),(0,\shortminus c_3),(0,c_3),(c_4,0)\} \;\;\forall \;\; \|b'\|>1 \;\;\; \lor \;\;\;  \{(\shortminus c_4,0),(\shortminus c_3,0),(c_3,0),(c_4,0)\}\;\; \forall\;\; \|b'\|<1$$
and has the following shape 
\newline

\newrgbcolor{qqzzff}{0 0.6 1}
\newrgbcolor{ttqqqq}{0.2 0 0}

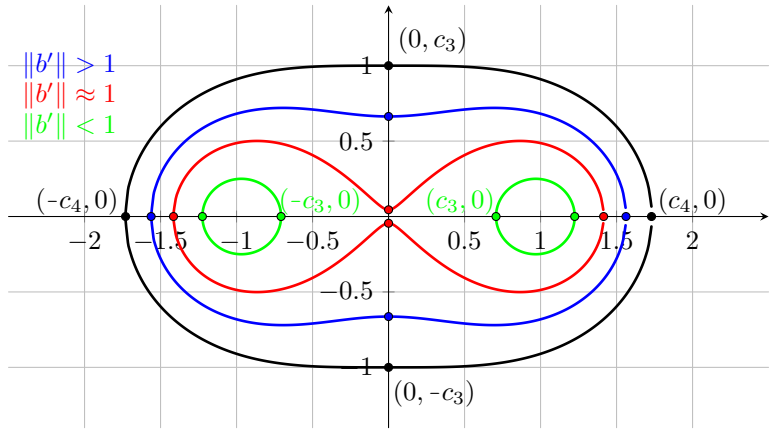
\begin{figure}[hbt!]
\centering
\begin{tikzpicture}[line cap=round,line join=round,>=triangle 45,x=2cm,y=2cm]
\begin{axis}[
x=2.0cm,y=2.0cm,
axis lines=middle,
ymajorgrids=true,
xmajorgrids=true,
xmin=-2.5,
xmax=2.5,
ymin=-1.4,
ymax=1.4,
xtick={-2,-3/2,...,2},
ytick={-1,-1/2,...,1},]
\clip(-2.6541760682099147,-6.42426466103129) rectangle (10.277881323330046,5.980780987279769);

\draw[line width=1pt,color=black,smooth,samples=200,domain=-1.7320492250169548:1.7320480088065726] plot(\x,{sqrt(-1-(\x)^(2)+sqrt(4*(\x)^(2)+(2/sqrt(2))^(4)))});
\draw[line width=1pt,color=black,smooth,samples=200,domain=-1.7320492250169548:1.7320480088065726] plot(\x,{-sqrt(-1-(\x)^(2)+sqrt(4*(\x)^(2)+(2/sqrt(2))^(4)))});

\draw[line width=1pt,color=blue,smooth,samples=200,domain=-1.5620482404410878:1.5620468802436978] plot(\x,{sqrt(-1-(\x)^(2)+sqrt(4*(\x)^(2)+(-1.2)^(4)))});
\draw[line width=1pt,color=blue,smooth,samples=200,domain=-1.5620482404410878:1.5620468802436978] plot(\x,{-sqrt(-1-(\x)^(2)+sqrt(4*(\x)^(2)+(-1.2)^(4)))});

\draw[line width=1pt,color=red,smooth,samples=200,domain=-1.4142102367228777:1.414212399303442] plot(\x,{sqrt(-1-(\x)^(2)+sqrt(4*(\x)^(2)+(1.001)^(4)))});
\draw[line width=1pt,color=red,smooth,samples=200,domain=-1.4142102367228777:1.414212399303442] plot(\x,{-sqrt(-1-(\x)^(2)+sqrt(4*(\x)^(2)+(1.001)^(4)))});

\draw[line width=1pt,color=green,smooth,samples=200,domain=0.7071079069455056:1.2247437735079298] plot(\x,{sqrt(-1-(\x)^(2)+sqrt(4*(\x)^(2)+(1.0001/sqrt(2))^(4)))});
\draw[line width=1pt,color=green,smooth,samples=200,domain=0.7071079069455056:1.2247437735079298] plot(\x,{-sqrt(-1-(\x)^(2)+sqrt(4*(\x)^(2)+(1.0001/sqrt(2))^(4)))});

\draw[line width=1pt,color=green,smooth,samples=200,domain=0.7071079069455056:1.2247437735079298] plot(-\x,{sqrt(-1-(\x)^(2)+sqrt(4*(\x)^(2)+(1.0001/sqrt(2))^(4)))});
\draw[line width=1pt,color=green,smooth,samples=200,domain=0.7071079069455056:1.2247437735079298] plot(-\x,{-sqrt(-1-(\x)^(2)+sqrt(4*(\x)^(2)+(1.0001/sqrt(2))^(4)))});

\begin{scriptsize}
\draw[color= green] (0.47,0.1) node {$(c_3,0)$};
\draw[color= green] (-0.45,0.1) node {$(\shortminus c_3,0)$};

\draw[color= black] (0.29,1.17) node {$(0,c_3)$};
\draw[color= black] (0.3,-1.17) node {$(0,\shortminus c_3)$};

\draw[color= black] (2,0.1) node {$(c_4,0)$};
\draw[color= black] (-2.05,0.1) node {$(\shortminus c_4,0)$};

\draw [fill= black] (0,1) circle (1.5pt);
\draw [fill= black] (0,-1) circle (1.5pt);
\draw [fill= black] (1.73,0) circle (1.5pt);
\draw [fill= black] (-1.73,0) circle (1.5pt);

\draw [fill= green] (0.7071,0) circle (1.5pt);
\draw [fill= green] (-0.7071,0) circle (1.5pt);
\draw [fill= green] (1.22474,0) circle (1.5pt);
\draw [fill= green] (-1.22474,0) circle (1.5pt);

\draw [fill= red] (1.41492,0) circle (1.5pt);
\draw [fill= red] (-1.41492,0) circle (1.5pt);
\draw [fill= red] (0,0.0447325) circle (1.5pt);
\draw [fill= red] (0,-0.0447325) circle (1.5pt);

\draw [fill= blue] (1.56205,0) circle (1.5pt);
\draw [fill= blue] (-1.56205,0) circle (1.5pt);
\draw [fill= blue] (0,0.663325) circle (1.5pt);
\draw [fill= blue] (0,-0.663325) circle (1.5pt);

\draw[color= blue]   (-2.1,1.0) node {$\|b'\|>1$};
\draw[color= red]    (-2.1,0.8)  node {$\|b'\| \approx 1$};
\draw[color= green] (-2.1,0.6) node {$\|b'\|<1$};

\end{scriptsize}

\end{axis}
\end{tikzpicture}
\caption{Cassini ovals for $a'=1$ with axis intersection points} 
\label{fig:F4}
\end{figure}

\subsection{Two\space $X$-axis intersections}
\label{section:CO1}
A relationship between the Cassini oval and the congruent number problem becomes apparent from a system of equations first defined by Kurt Heegner \cite{8} in $1952$ on page $228$. The key point is to see that his variables $c_1,c_2$ can be defined in function of the variables $N,f_1,f_2$ from section  $\ref{section:Ellipse}$ and that the Cassini oval intersection points are defined by his variables $c_3,c_4$ from the equations 
\begin{eqnarray} 
 (c_1,c_2) = & \left(\|f_1 f_2\| , \frac{1}{2}\|N f_1^2 -f_2^2\|\right) \\
 (c_3,c_4) = &\left(\sqrt{\|N c_1^2-c_2^2\|}, \sqrt{N c_1^2+c_2^2} \right)
\end{eqnarray} 
The right triangle for the congruent number $N$  is then represented by
\begin{equation} 
(a,b,c)=\left(\frac{c_3 c_4}{c_1 c_2},\frac{2 c_1 c_2 N}{c_3 c_4},\sqrt{\frac{c_3^2  c_4^2}{c_1^2 c_2^2}+\frac{4 c_1^2  c_2^2 N^2}{c_3^2 c_4^2}}\right)
\end{equation} 
The corresponding Cassini oval is obtained by assigning $(a',b')=\left(c_2,c_1\sqrt{N}\right)$ resulting in the equation
\begin{equation}
y(x)=\sqrt{\sqrt{4c_2^2 x^2+c_1^4 N^2 } - x^2 - c_2^2}
\end{equation}
For example for $N=29$ with $(f_1,f_2)=(1,\shortminus 13)$ we obtain $(c_1,c_2,c_3,c_4)=(13,70,1,99)$
for the right triangle 
$$(a,b,c)=\left(\frac{99}{910}, \frac{52780}{99}, \frac{48029801}{90090}\right)$$ 
and the Cassini oval 
$$ y(x)=\sqrt{\sqrt{4*70^2 x^2+13^4* 29^2 } - x^2 - 70^2} $$ 
which is almost a lemniscate with the intersection points 
$$\{(\shortminus 99,0),(0,\shortminus1),(0,1),(99,0)\}$$
The plots for the examples are added in the Sagemath notebook $\texttt{05\_Cassini.ipynb}$.
\pagebreak
\newline
Similarly to section $\ref{section:Ellipse}$ we can adjoin $f_2$ for certain congruent numbers.
\newline
\linebreak
For $N=79$ with $(f_1,f_2)=\left(125, 52\right)$ adjoining $f_2$ by $\sqrt{N}$  we get
$$(c_1,c_2,c_3,c_4)=\left(6500\sqrt{79}, \frac{1020759}{2}, \frac{12719}{2}\sqrt{79}, \frac{1447991}{2}\right)$$
Even though $c_1$ and $c_3$ are not rational we obtain the rational right triangle 
$$(a,b,c)=\left(\frac{233126551}{167973000}, \frac{335946000}{2950969}, \frac{56434050774922081}{495683115837000} \right)$$
due to the cancellation of the squareroots in $(65)$ and the corresponding Cassini oval is 
$$y(x)=\sqrt{\sqrt{1020759^2 x^2+6500^4 79^4 } - x^2 - \frac{1020759^2}{4}}$$
By using a second system of equations a Cassini oval with four $X$-axis intersections for $N=79$
\newline
such that $(c_1,c_2,c_3,c_4)$ are integers is shown in the next section. 
\newline
\linebreak
For $N=62$  with $(f_1,f_2)=\left(20, 7\right)$ adjoining $f_2$ by $\sqrt{2N}$  we get
$$(c_1,c_2,c_3,c_4)=\left(280\sqrt{31}, 9362, 1426\sqrt{31}, 15438\right)$$
the rational triangle 
$$(a,b,c)=\left(\frac{177537}{21140}, \frac{84560}{5727}, \frac{2056525601}{121068780} \right)$$
and the corresponding Cassini oval
$$y(x)=\sqrt{\sqrt{4*9362^2 x^2+280^4 31^2 61^2 } - x^2 - 9362^2}$$
\subsection{Four $X$-axis intersections}
Another system of equations relating the congruent number problem to the Cassini oval is 
\begin{eqnarray} 
 (c_3,c_4) = &\left(f_1^2-f_2^2, 2 f_1 f_2\right) \\
 (c_1,c_2) = & \left(\sqrt{\frac{c_4^2-c_3^2}{N}}, f_1^2 +f_2^2\right) 
\end{eqnarray} 
where the oval 
\begin{equation}
y(x)=\sqrt{\sqrt{8 c_2^2 x^2 + c_1^4 N^2} - 2 x^2 - c_2^2}
\end{equation}
splits into two separate closed loops around a focus with the intersection points 
$$ \{(-c_4,0),(-c_3,0),(c_3,0),(c_4,0)\}$$  
representing the right triangle 
\begin{equation}
(a,b,c)=\left(\frac{c_1 c_2 N}{c_3 c_4},\frac{2 c_3 c_4}{c_1 c_2},\sqrt{\frac{c_1^2 c_2^2 N^2}{c_3^2 c_4^2}+\frac{4 c_3^2 c_4^2}{c_1^2 c_2^2}}\right)
\end{equation}
For the example $N=79$  with $(f_1,f_2)=\left(125, 52\right)$ we now obtain the Cassini oval 
$$y(x)=\sqrt{\sqrt{8 * 18329^2 x^2 + 161^4 79^2} - 2 x^2 - 18329^2}$$
with the four $X$-axis intersection points  
$$\{ (\shortminus13000, 0), (\shortminus12921, 0) , (12921, 0),(13000, 0)\}$$
\newline
And for the example $N=62$ with $(f_1,f_2)=\left(20, 7 \sqrt{2}\right)$ we obtain the Cassini oval 
$$y(x)=\sqrt{\sqrt{8 * 498^2 x^2 + 4*23^4 62^2} - 2 x^2 - 498^2}$$
with  the intersection points  
$$\{ (\shortminus 280 \sqrt{2}, 0), (\shortminus302, 0) , (302, 0),(280 \sqrt{2}, 0)\}$$
\pagebreak
\section{The tangent method}
The Tangent method generates a solution set $S_n=\{(f_1,f_2)_1,..,(f_1,f_2)_n\}$ for a congruent number $N$  from an initial known right triangle $(a,b,c)$ and is similar to the elliptic curve point doubling method.
\newline
\linebreak
The solutions $S_i$ correspond to the intersections of the tangent lines of the elliptic curve  $E_N: y^2=x^3-N^2x$ such that the next tangent line touches $E_N$ at the intersection point found by the previous tangent line.
\newline
\linebreak
From the known hypotenuse $c$ of a right triangle for a congruent number $N$ we assign $(c_1,c_2)$ such that 
$$(c_1,c_2)=\left( Denominator[c],\frac{1}{2}Numerator[c]\right)$$ 
for which we obtain a new solution $S_i=(f_1,f_2)_i$ by solving the binary quadratic form for $c_2$ from the equations
\begin{equation} 
 (c_1,c_2) = \left(\|f_1 f_2\| , \frac{1}{2}\|N f_1^2 -f_2^2\|\right) 
\end{equation} 
The solution $S_i$ corresponds to a new right triangle $(a,b,c)_i$ by evaluating the equations defined in section  $\ref{section:CO1}$
Repeating this method we obtain different right triangles for the same congruent number $N$.
\newline
\linebreak
For example for $N=5$ starting from the hypotenuse $c=\frac{41}{6}$ we obtain the following solutions
\begin{align*}  
S_{1-4}=\{(&3,2) ,(372,2009),(169317668184,15811196552161),(1336220772668316930638357029463135419039\backslash\\ & 997035301712,62496947695267799013412096545625364258488963961427841)\} 
\end{align*}
The figure below illustrates the method showing the first two tangent lines for the elliptic curve $E_5: y^2=x^3-25 x$  and the rational points $P_{1,2,3}$.
The first tangent line, touches $E_5$ at the point $P_1$ and intersects $E_5$ at the point $P_2$. Using $P_2$ as the next tangent point the rational point $P_3$ is obtained.
\newline
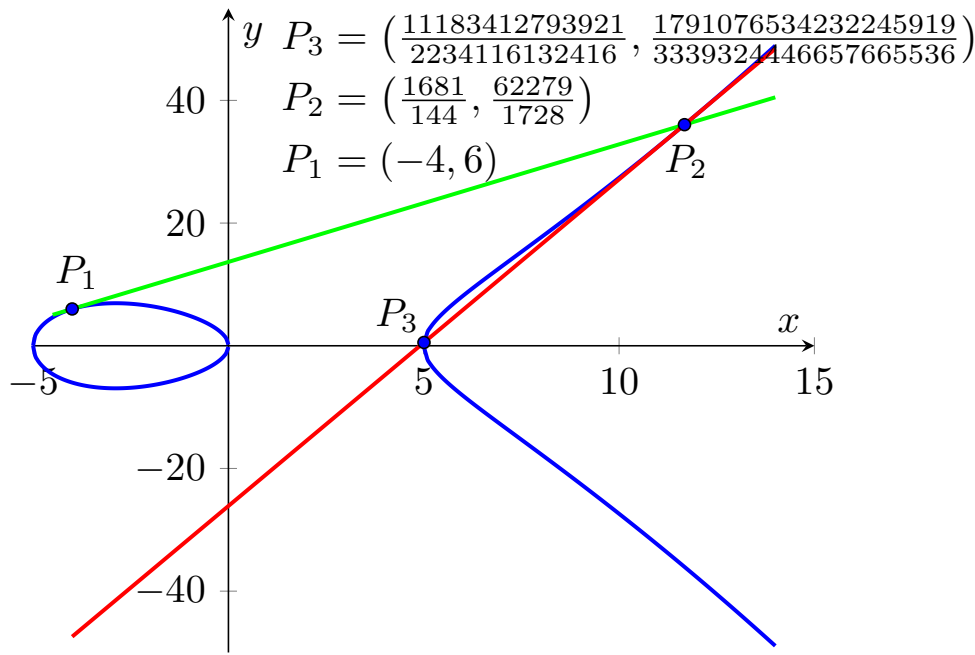
\begin{figure}[hbt!]
\centering
\begin{tikzpicture}[scale=1.5]
\begin{axis}[
        xmin=-5,
        xmax=15,
        ymin=-50,
        ymax=55,
        xlabel={$x$},
        ylabel={$y$},
        axis lines=middle,
        samples=201,
        smooth,   
        clip=false,
    ]
   
\draw[line width=1pt,color=blue,smooth,samples=100,domain=5:14] plot(\x,{sqrt((\x)^(3)-5^(2)*(\x)});
\draw[line width=1pt,color=blue,smooth,samples=100,domain=-5:0] plot(\x,{sqrt((\x)^(3)-5^(2)*(\x)});
\draw[line width=1pt,color=blue,smooth,samples=100,domain=5:14] plot(\x,{-sqrt((\x)^(3)-5^(2)*(\x)});
\draw[line width=1pt,color=blue,smooth,samples=100,domain=-5:0] plot(\x,{-sqrt((\x)^(3)-5^(2)*(\x)});

\draw[line width=1pt,color=green,smooth,samples=100,domain=-4.5:14] plot(\x,{1.9167*(\x)+13.667});
\draw[line width=1pt,color=red,smooth,samples=100,domain=-4:14] plot(\x,{5.3248*(\x)-26.118});

\begin{scriptsize}

\draw[color= black]   (1,30) node[anchor = west] {$P_1=(-4,6)$};
\draw[color= black]   (1,40) node[anchor = west] {$P_2=\left(\frac{1681}{144},\frac{62279}{1728}\right)$};
\draw[color= black]   (1,50) node[anchor = west] {$P_3=\left(\frac{11183412793921}{2234116132416},\frac{1791076534232245919}{3339324446657665536}\right)$};

\draw[color= black]   (-4.8,12) node[anchor = west] {$P_1$};
\draw[color= black]   (10.8,30) node[anchor = west] {$P_2$};
\draw[color= black]   (3.4,5) node[anchor = west] {$P_3$};

\draw [fill= blue] (-4,6) circle (1.5pt);
\draw [fill= blue] (11.674,36.041) circle (1.5pt);
\draw [fill= blue] (5.0057,0.53636) circle (1.5pt);
\end{scriptsize}
\end{axis}
\end{tikzpicture}
\caption{The tangent method} 
\label{fig:F5}
\end{figure}
\newline
It should be noted that the solutions grow quite fast like for example for $N=79$ the first three solutions are
\begin{align*}  
S_{1-3}=\{(&2080281,238277000), (55260645511189879706636193594000,3223389202505003051748398476629439), \\ 
(&343513126900812678740983883808614992645462095396485294670725007133093011002487199141113\backslash\\ &  5222198071436263742039481934242611258578572000 \\
,&107898845343126533895905583080797449600586402309215336477362300480422977554420331904855\backslash\\ & 288440002098386935913690366578027064421253187841)\}
\end{align*}
The equations and examples can be verified in the Sagemath notebook $\texttt{06\_Tangent\_method.ipynb}$.
\pagebreak 
\section{The prime footprints }
\label{section:FP}
For certain congruent numbers $N$ which are prime or 2 times a prime,  one can define footprint equations for the variables $p,q$ in function of $N,m,n \in \mathbb{Z}$ for the right triangle with sides $a,b,c \in \mathbb{Q}$ 
\begin{equation}(a,b,c)=\left(\frac{p}{q},2 N\frac{q}{p},\frac{\sqrt{p^4+4N^2q^4}}{pq}\right)\end{equation}
The equations depend on the congruence $mod \;8$ and allow to see a relationship to a combination of conics.
\newline
Their solutions $m,n$ for primes $P<1000$, are listed in the Tables $\ref{section:TFP}$
\newline 
\linebreak
Table $\ref{table:T0}$ for $N = P \;\; \not \forall \;\; P \equiv 1 \mod  8$  
\newline
\indent * For the primes of this kind represented by $ N=(m^4 + 6 m^2 n^2 + n^4)$ we have   
\begin{equation} 
(p, q) =\left(N(m^2 - n^2), 2 m n (m^2+ n^2) \right) 
\end{equation}
\indent * For the other congruent primes of this kind we have 
\begin{equation} 
(p, q) =\left(\frac{1}{4}(m^2 + n^2) \sqrt{N\left(\left(2m n\right)^2-\left(m^2-n^2\right)^2\right)}, \frac{1}{2} m n (m^2 - n^2)\right) 
\end{equation}
Table $\ref{table:TI}$ for $N = P \;\; \forall \;\; P \equiv 5 \mod  8$ 
\begin{equation} 
(p, q) =\left(\sqrt{(m^2n^2N)^2 - \frac{(m^2N-n^2)^4}{16}}, \frac{mn(m^2N-n^2)}{2}\right) 
\end{equation}
Table $\ref{table:TII}$ for $N = P \;\; \forall\;\; P \equiv 7 \mod  8$ 
\begin{equation} (p, q) = \left((m^2+n^2)\sqrt{N\left(\left(2 m n\right)^2-\left(m^2-n^2\right)^2\right)},2mn\left(m^2-n^2\right)\right)\end{equation}
Table $\ref{table:TIII}$ for $N = 2P\;\;  \forall\;\; P \equiv 7 \mod  8$ 
\begin{equation}(p, q) = \left((m^2+2n^2)\sqrt{N\left(\left(2 \sqrt{2} m n\right)^2-\left(m^2-2n^2\right)^2\right)},2\sqrt{2}mn\left(m^2-n^2\right)\right)
\end{equation}
Table $\ref{table:TIV}$  for  $N =2P\;\;  \forall\;\;  P \equiv 3 \mod 8 $ 
\begin{equation} 
(p, q) = \left(\left(m^2-n^2-2mn\right)\sqrt{\frac{N}{2} \left(\left(m-n\right)^2+2 m^2\right)\left(\left(m+n\right)^2 + 2 n^2\right)},\left(m^2-n^2+2mn\right)\left(m^2+n^2\right)\right) 
\end{equation}
Evaluating one example from each table we get the following right triangles :
\begin{align*}
  (N,m,&n) =(353,4,1)  \\ 
(a,b ,&c) = \left(\frac{5295}{136}, \frac{272}{15}, \frac{87617}{2040}\right) \\
  (N,m,&n) =(761,31,51)  \\ 
(a,b ,&c) = \left(\frac{66411709}{1296420},  \frac{2592840}{87269},
  \frac{6699926952721}{113137276980}\right) \\ 
  (N,m,&n) =(173, 10865, -343141)  \\ 
(a,b ,&c) = \left(\frac{418416739097462232963}{181421867613059954270},  \frac{62771966194118744177420}{418416739097462232963},
  \frac{11389552969201600543101928087171460571651881}{75909946247628040203029119534348866602010}\right) \\ 
    (N,m,&n) =(191,  27469, 11580)  \\ 
 (a,b ,&c) =  \left(\frac{1726816796630813713}{394718867434084440},  \frac{789437734868168880}{9040925636810543},
 \frac{311996818759910472998178689881743841}{3568623927917636168751328944250920} \right) \\ 
   (N,m,&n) =(382, 540, 239)  \\ 
 (a,b ,&c) =  \left(\frac{447382566673}{11444911740},  \frac{45779646960}{2342317103},
 \frac{1171595729834345971681}{26807612510927489220}\right) \\ 
  (N,m,&n)  =(326 , 170, -69)  \\ 
 (a,b ,&c) =  \left(\frac{28931957373}{22855819},   \frac{91423276}{177496671},
 \frac{5135326544339012645}{4056831785478549}\right)
\end{align*}
The equations tables and examples can be verified in the Sagemath notebook $\texttt{07\_Prime\_footprints.ipynb}$.
\pagebreak
\section{The unseen recurrence}
A general recurrence hidden in the equation of Pythagoras defines infinite trees of congruent number sequences and their corresponding Pythagorean triples.
\newline
\linebreak
In order to see this we use the right triangle with sides $a,b,c$ for a congruent number $N$ which can be defined in terms of $p,q \in \mathbb{Z}$ by 
\begin{equation}(a,b,c)=\left(\frac{p}{q},2 N\frac{q}{p},\frac{\sqrt{p^4+4N^2q^4}}{pq}\right)\end{equation}
The congruent number recurrence is then defined by  
\begin{equation}(N, p, q) = \left(\sqrt{p^4 + 4 N^2 q^4}, p  \sqrt{p^4 + 4 N^2 q^4}  ,  q^2 N\right)\end{equation}
Evaluating the right hand side using an initial solution $(a,b,c)$ defined by $N,p,q$ we obtain a new set $N,p,q$ resulting in a new congruent number with its corresponding Pythagorean triple. By this iteration one generates infinite sequences of congruent numbers.
It is important to notice that the sequences and triples change depending on the side $a$ or $b$ one chooses to define the variables $N,p,q$ not only initially but also during the iteration process resulting in a tree structure.
\newline
\linebreak
So by using any triple definition $(a,b,c)$ and assigning $(N,p,q)$ to 
$$ \begin{array}{cccc}& (N,p,q)_a &= & (N,Numerator[a],Denominator[a]) 
\\  or & & &
\\ & (N,p,q)_b & = & (N,Numerator[b],Denominator[b])
\\ \end{array} $$
the recurrence defines infinite tree-sets of Pythagorean triples for which $N$ is a congruent number sequence 
by choosing a tree-side-walk iteration, denoting the chosen side-path sequence in the subscript of the triple $(a,b,c)$. 
\newline
\linebreak
For example applying the recurrence to Euclid's definition we get for the side-paths $a^i,ab,bb$ the following triples
$$(A,B,C)_{a^i}=\left(\frac{m^{2^{i+1}}-n^{2^{i+1}}}{(m n)^{2^{i-1}}},2 (m n)^{2^{i-1}},\frac{m^{2^{i+1}}+n^{2^{i+1}}}{(m n)^{2^{i-1}}}\right)$$
$$(A,B,C)_{ab}=\left(\frac{4 m n \left(m^2+n^2\right)}{m^2-n^2},m^2-n^2,\frac{m^4+6 m^2 n^2+n^4}{m^2-n^2}\right)$$ 
$$(A,B,C)_{bb}=\left(\frac{8 m n \left(m^6+7 m^4 n^2+7 m^2 n^4+n^6\right)}{\left(m^2-n^2\right)^2},\left(m^2-n^2\right)^2,\frac{m^8+28 m^6 n^2+70 m^4 n^4+28 m^2 n^6+n^8}{\left(m^2-n^2\right)^2}\right)$$
Another way to illustrate the recurrence is to write the congruent number followed by the triple by $N:(a,b,c)$ and apply the recurrence denoting the choice of the side-walk iteration by $\overset{a}{\rightarrow}$ or $ \overset{b}{\rightarrow}$ evaluating the next congruent number and its corresponding triple.
\newline
\linebreak
For the right triangles from our first example for $(m,n)=(2,1)$ the first few iterations gives us 
$$ \begin{array}{ccccccccccc}
6&:&(3, 4, 5)
     & \overset{a}{\rightarrow}&    15&:& \left(\frac{15}{2},4,\frac{17}{2}\right)  & \overset{a}{\rightarrow}&    {255}&:&\left(\frac{255}{4},8,\frac{257}{4}\right)
\\&  &  &  &  & &  & \overset{b}{\rightarrow} &   {34}&:&\left(\frac{136}{15},\frac{15}{2},\frac{353}{30}\right)
\\&  &  & \overset{b}{\rightarrow}&    20&:& \left(\frac{40}{3},3,\frac{41}{3}\right) & \overset{a}{\rightarrow} &   {1640}&:&\left(\frac{3280}{9},9,\frac{3281}{9}\right)
\\&  &  &  &  & &  & \overset{b}{\rightarrow} &   {41}&:&\left(\frac{123}{20},\frac{40}{3},\frac{881}{60}\right)

\\ 34&:&(\frac{15}{2}, \frac{136}{15} \frac{353}{30})   
     & \overset{a}{\rightarrow}&     353&:& \left(\frac{5295}{136},\frac{272}{15},\frac{87617}{2040}\right)  & \overset{a}{\rightarrow}&   {30928801}&:&\left(\frac{463932015}{18496},\frac{36992}{15},\frac{6992534657}{277440 }\right) 
\\&  &  &  &  & &  & \overset{b}{\rightarrow} &   175234&:&\left(\frac{47663648}{79425},\frac{79425}{136},\frac{9045146753}{10801800}\right)
\\ & &  & \overset{b}{\rightarrow} &   24004&:& \left(\frac{96016}{225},\frac{225}{2},\frac{198593}{450}\right)  & \overset{a}{\rightarrow}&   9534052744&:&\left(\frac{38136210976}{50625},\frac{50625}{2},\frac{76315468673}{1012
   50}\right)
\\&  &  &  &  & &  & \overset{b}{\rightarrow} &   198593&:&\left(\frac{44683425}{96016},\frac{192032}{225},\frac{21001035137}{21603600}\right)

\\ 41&:&(\frac{40}{3}, \frac{123}{20}, \frac{881}{60}) 
     & \overset{a}{\rightarrow} &  1762&:& \left(\frac{70480}{369},\frac{369}{20},\frac{1416161}{7380}\right)  & \overset{a}{\rightarrow}&   4990551364&:&\left(\frac{199622054560}{136161},\frac{136161}{20},\frac{3992484137921
   }{2723220}\right)
\\&  &  &  &  & &  & \overset{b}{\rightarrow} &   1416161&:&\left(\frac{522563409}{704800},\frac{1409600}{369},\frac{1012025897921}{260071200}\right)
\\ &  &  &\overset{b}{\rightarrow} &  36121&:& \left(\frac{108363}{400},\frac{800}{3},\frac{456161}{1200}\right) &  \overset{a}{\rightarrow}&  16476991481&:&\left(\frac{49430974443}{160000},\frac{320000}{3},\frac{156882857921}{4
   80000}\right)
\\&  &  &  &  & &  & \overset{b}{\rightarrow} &   912322&:&\left(\frac{729857600}{325089},\frac{325089}{400},\frac{310482857921}{130035600}\right)

\\ 7&:&(\frac{24}{5}, \frac{35}{12}, \frac{337}{60})  
     & \overset{a}{\rightarrow} &   674&:& \left(\frac{16176}{175},\frac{175}{12},\frac{196513}{2100}\right) &  \overset{a}{\rightarrow}&   264899524&:&\left(\frac{6357588576}{30625},\frac{30625}{12},\frac{76296827713}{3675
   00}\right)
\\&  &  &  &  & &  & \overset{b}{\rightarrow} &   196513&:&\left(\frac{34389775}{97056},\frac{194112}{175},\frac{19777624897}{16984800}\right)
\\ &  &  & \overset{b}{\rightarrow} &  2359&:& \left(\frac{11795}{144},\frac{288}{5},\frac{72097}{720}\right) 
 &  \overset{a}{\rightarrow}&  170076823&:&\left(\frac{850384115}{20736},\frac{41472}{5},\frac{4338014017}{103680}
   \right)
\\&  &  &  &  & &  & \overset{b}{\rightarrow} &   144194&:&\left(\frac{41527872}{58975},\frac{58975}{144},\frac{6917904193}{8492400}\right)
\\ 
\end{array} $$
The equations and examples can be verified in the Sagemath notebook $\texttt{08\_Unseen\_recurrence.ipynb}$.
\pagebreak
\section{The  Fibonacci/Lucas numbers}
We define two congruent number sequences related to the Fibonacci/Lucas numbers and show how the congruent number solution for $5$ first found by Fibonacci can be expressed in terms of the Fibonacci/Lucas numbers. 
\newline
\linebreak
Combining the identitiy between the Fibonacci \cite{9} and the Lucas numbers \cite{10}
\begin{equation} L_{n}^2= 5 F_{n}^2+4(-1)^n\end{equation}
together with the algebraic equation
\begin{equation} \left(L_n^2-4\right)^2+\left(4L_n\right)^2= \left(L_n^2+4\right)^2\end{equation}
we can define one right triangle for $n$ even and one for $n$ odd due to the $(-1)^n$ in the identity.
\subsection{The even indices}
The right triangle resulting from the algebraic equation and the identity for even indices is 
\begin{equation}(a,b,c)_n'=(5 F_{2n}^2, 4L_{2n},L_{2n}^2+4) \end{equation}
Reducing the sides by $F_{2n}$ we obtain the rational right triangle  
\begin{equation}(a,b,c)_n=\left(5 F_{2n}, 4\frac{L_{2n}}{F_{2n}},\frac{L_{2n}^2+4}{F_{2n}}\right)  \end{equation}
and  the area for this triangle defines the congruent number sequence
\begin{equation} N_{n}=10 L_{2n}= 10 \left(F_{2n-1}+F_{2n-1}\right)\end{equation}
An interesting historical example is $N_3 = 10 L_6 = 5 * 6^2$ ,  for which the sides can be reduced by $6$ due to the square factor obtaining a right triangle with area $5$ in terms of the Fibonacci/Lucas numbers
\begin{equation}(a,b,c)=\left(\frac{5 F_6}{6},\frac{4 L_6}{6 F_6},\frac{L_6^2+4}{6 F_6}\right)=\left(\frac{20}{3},\frac{3}{2},\frac{41}{6}\right)\end{equation} 
Multiplying the identity (81) for $n$ even by $2000$ we obtain   
$$   20 (10 L_{2n})^2 = (10^2 F_{2n})^2 + 20^3 $$ 
from which follows that there exists a parallel line to the $y$-axis at $x=-20$ intersecting and connecting all the
 Elliptic curves $E_{N_{n}}: y^2=x^3- N_{n}^2 \;x$ for this congruent number sequence 
at the rational points   
\begin{equation}
P_0= \left(-20, 10^2 F_{2n}\right) 
\end{equation}
These points are closely related to the two points 
\begin{equation}
P_1= \left(\frac{1}{2}a(a + c) , \frac{1}{2}a^2(a + c)\right) \;\;\; ,\;\;\; P_2 = \left(\frac{c^2}{4}, \frac{c}{8} (a^2 - b^2)\right) 
\end{equation}
By adding the trivial point $(0,0)$  to $P_0$ we obtain the point $P_1$ and the point $P_2$ is obtained by doubling point $P_0$
\begin{equation}
P_1= (0,0) + P_0  \;\;\;,\;\;\; P_2 = 2 P_0
 \end{equation}
\subsection{The odd indices}
The right triangle resulting from the algebraic equation and the identity for odd indices is 
\begin{equation}(a,b,c)_n=(L_{2n+1}^2-4,4L_{2n+1},5 F_{2n+1}^2) \end{equation}
giving the congruent number sequence
\begin{equation}N_{n}=2\left(L_{2n+1}^2-4\right) L_{2n+1} \end{equation}
\newline
\linebreak
The Elliptic curve  $E_{N_{n}}: y^2=x^3- N_{n}^2 \;x$ for this congruent number sequence 
also has the rational points of infinite order as defined by $(88)$.
\newline
\linebreak
The equations and examples can be verified in the Sagemath notebook $\texttt{09\_Fibonacci\_Lucas.ipynb}$.
\pagebreak
\section{The Chebyshev polynomials}
We define a congruent number family related to the Chebyshev polynomials and show their relationship to the Heron-Brahmagupta triangles, their semiperimeter and some properties of their associated elliptic curves.
\newline
\linebreak
Using the Pell identitiy \cite{11} relating the Chebyshev polynomials of the first $T_{m}(n)$ and second kind $U_{m-1}(n)$
\begin{equation} T_{m}(n)^2= (n^2-1) U_{m-1}(n)^2+1\end{equation}
in combination with the algebraic equation for the polynomials of the first kind
\begin{equation} \left(T_{m}(n)^2-1\right)^2+ 4 T_{m}(n)^2 =  \left(T_{m}(n)^2+1\right)^2\end{equation}
we can define the right triangle $(a,b,c)'$ for   \(m>0\) and \(n>1\)  by 
\begin{equation}(a,b,c)_{m,n}'=\left(\left(n^2-1\right) U_{m-1}\left(n\right)^2, 2 T_{m}(n),T_{m}\left(n\right)^2+1\right) \end{equation}
Reducing the sides due to the square factor of $U_{m-1}(n)$ we get the rational triangle
\begin{equation}(a,b,c)_{m,n}=\left(\left(n^2-1\right) U_{m-1}(n),\frac{2 T_m(n)}{U_{m-1}(n)} ,\frac{T_m(n)^2+1}{U_{m-1}(n)} \right)\end{equation}
and the area $N_{m,n}$ for this triangle defines the congruent number polynomial
\begin{equation}N_{m,n}= (n^2-1) T_{m}(n)\end{equation}
related to the Chebyshev polynomial of the first kind. 
\newline
The corresponding congruent number Elliptic curve  $E_{N_{m,n}}: y^2=x^3- N_{m,n}^2 \;x$ has the rational point 
\begin{equation}
P_0:\left(1-n^2, (n^2-1)^2U_{m-1}(n)\right)
\end{equation}
and also the two rational points $P_1,P_2$ as defined by $(88)$ with the same relationships as $(89)$.
\subsection{The Brahmagupta triangles}
The Brahmagupta triangles are Heronian triangles with consecutive integer sides $A,B,C$ ,  integer semiperimeter $P$ and integer area $S$.
A generalization of these triangles as defined in \cite{12} with consecutive integer sides is 
\begin{eqnarray}(A,B,C)_{t} & = &\left(  t -1 ,  t  ,  t + 1\right)
\\ P_{t} & = &\frac{1}{2}(A+B+C)=3\frac{t}{2}
\\ S_{t} &= &\sqrt{P(P-A)(P-B)(P-C)}=  \frac{t}{2} \sqrt{3 \left(\left(\frac{t}{2}\right)^2-1\right)} 
\end{eqnarray}
The Brahmaputa triangles are obtained replacing $t$ by $2T(t,2)$ where $T(t,2)$ is the Chebyshev polynomials of the first kind. For these triangles the semiperimeter $P_t$ and the area $N_{m,n}$ defined by $(96)$ for a right triangle $(95)$ are related by 
\begin{equation}
P_{t}=N_{t,2} 
\end{equation}
The elliptic curves for the generalized triangles  $$E: y^2 = (x + A B) (x + B C) (x + A C) $$ 
have a minimum rank of $1$ due to the integral points $Q_{0,1,2,3}:(x,y)$ of infinte order for $t>2$ 
\begin{equation}
Q_0:  \left(0,A B C\right) \;\;, \;\;
Q_1:  \left(-B^2, B \right) \;\;,\;\; 
Q_2:  \left(2 - A B , 2 C \right) \;\;,\;\;
Q_3:  \left(2 - B C , 2 A \right)  
\end{equation}
For $t=2$  the area of the triangle is $0$ because of the sides $(A,B,C)=(1,2,3)$. In this case the Elliptic curve reduces to $E : y^2= x^3+11 x^2+36 x +36$ which has rank $0$ and the points $Q_{0..3}$ have order $4$.
\newline
\linebreak
For example for $t=2 T(3,2)$ we obtain the Brahmagupta triangle with integer sides area and semiperimeter 
$$(A,B,C)=(51,52,53)  \;\;\;,\;\;\;  S=1170 \;\;\;,\;\;\; P=78$$ 
and the integral points $Q_{0,1,2,3}$ on the elliptic curve $E: y^2=(x+2652)(x+2756)(x+2703)$ are
$$ Q_0:(0,140556) \;\;,\;\; Q_1:(-2704,52) \;\;,\;\; Q_2:(-2650,106) \;\;,\;\; Q_3:(-2754,102)$$ 
The semiperimeter is a congruent number due to $(101)$ representing the area of the triangle defined by $(95)$ 
$$(a,b,c)=\left(45,\frac{52}{15},\frac{677}{15}\right) $$
with the three rational points $P_{0,1,2}$ on the elliptic curve  $E_{N_{3,2}}: y^2=x^3 - 78^2x$
$$ P_0:(-3,135) \;\;,\;\; P_1:(2028,91260) \;\;,\;\; P_2:\left(\frac{458329}{900},\frac{306627517}{27000}\right) $$ 
The equations and examples can be verified in the Sagemath notebook $\texttt{10\_Chebyshev\_polynomials.ipynb}$.
\pagebreak
\section{The Fermat  triangles}
A naive recursive method is presented to generate solutions to a problem posed by Fermat in $1643$ finding the smallest right triangle with sides $$(a,b,c)=(4565486027761,   1061652293520 ,  4687298610289)$$ for which the sum of the sides $a+b = 23721592^2 $ and the hypotenuse $c = 21650172^2$ are squares. 
\newline
\linebreak
The Pythagorean triple for this solution is generated using Euclid's formula assigning $(m,n) = (2150905, 246792)$.
 \newline
 \linebreak
The method is based on the fact that the triple solutions for which the hypotenuse $c$ and the sum $a + b$ or the difference $a - b$ are squares form an infinite tree.
\newline
\linebreak
Denoting the solutions for which the sum is a square by $P_i$ and the solutions for which the difference is a square by $N_i$ the following table shows the solutions $(a,b,c)$ for the first and second $N_i$ and $P_i$  triangles.
\newline
\begin{center}
\renewcommand{\arraystretch}{1.2}
\begin{tabu}{|c|c|c|c|c|c|c|c|c|c|} \hline
\multicolumn{2}{|c|[1pt]}{$$}& \multicolumn{2}{|c|[1pt]}{$N_1$}  & \multicolumn{2}{|c|[1pt]}{$N_2$} & \multicolumn{2}{|c|[1pt]}{$P_1$}  & \multicolumn{2}{|c|[1pt]}{$P_2$} \\ \hline \tabucline[1pt]\hline
\multicolumn{2}{|c|[1pt]}{$a$} & \multicolumn{2}{|c|[1pt]}{$-119$} & \multicolumn{2}{|c|[1pt]}{$2276953$} & \multicolumn{2}{|c|[1pt]}{$4565486027761$} & \multicolumn{2}{|c|[1pt]}{$214038981475081188634947041892245670988588201$} \\ \hline 
\multicolumn{2}{|c|[1pt]}{$b$} & \multicolumn{2}{|c|[1pt]}{$120$} & \multicolumn{2}{|c|[1pt]}{$-473304$} & \multicolumn{2}{|c|[1pt]}{$1061652293520$} & \multicolumn{2}{|c|[1pt]}{$109945628264924023237017010068507003594693720$} \\  \hline
\multicolumn{2}{|c|[1pt]}{$c$}  & \multicolumn{2}{|c|[1pt]}{$169$} & \multicolumn{2}{|c|[1pt]}{$2325625$}& \multicolumn{2}{|c|[1pt]}{$4687298610289$} & \multicolumn{2}{|c|[1pt]}{$240625698472667313160415295005368384723483849$} \\ \hline
\end{tabu}
\newline
\renewcommand{\arraystretch}{1}
\end{center}
The infinite tree starts from a root node where each node defines 2 child solutions.
\newline
\linebreak
The root node is defined by the fraction $x=\frac{p}{q}=\frac{1}{1}$ and represents the initial trivial solution $1^2+0^2 =1^2$.
\newline
\linebreak
Assigning $(p,q)=\left(Numerator[x],Denominator[x]\right)$ we can define the Pythagorean triple 
\begin{equation}(a, b, c) = \left(p\; q , -\frac{1}{2}(p^2 - q^2),\frac{1}{2} (p^2 + q^2)\right)\end{equation}
and determine the intermediate variables $(n,m)$
\begin{equation}(n, m) = (a - b , \sqrt{a + b}\sqrt{c})\end{equation}
\newline
These allow us to evaluate the fractions $(x_1,x_2)$ 
\begin{equation}(x_1, x_2) = \left(\frac{(2 m n)^2 + n^4 + 4 m n \sqrt{8 m^4 + n^4}}{ 16 m^4 + n^4}, \frac{(2 m n)^2 + n^4 - 4 m n \sqrt{8 m^4 + n^4}}{ 16 m^4 + n^4}\right)\end{equation}
\newline
generating the two child solutions by assigning $p,q$ to the numerator and denominator of the respective fractions for two new triples in case $p,q \ne 1$.
\newline
\linebreak
Inso this method steps from one solution to the next without any search.
\subsection{The recursive algorithm}
The method can be implemented into a recursive algorithm as follows. 
\begin{algorithm}
\caption{Fermat recursive algorithm}
\textbf{procedure} $Fermat(x)$ \tikzmark{right}\Comment{ $x \in \mathbb{Q}$  }
\begin{algorithmic} 
\State $(p,q) \leftarrow(Numerator(x),Denominator(x))$
\State $(a, b, c)  \leftarrow  (p\; q , -\frac{1}{2}(p^2 - q^2),\frac{1}{2} (p^2 + q^2))$  \tikzmark{right}\Comment{$Pythagorean$ triple }
\State $(n, m)  \leftarrow (a - b , \sqrt{a + b}\sqrt{c})$
\State $(x_1, x_2)   \leftarrow \left(\frac{(2 m n)^2 + n^4 + 4 m n \sqrt{8 m^4 + n^4}}{ 16 m^4 + n^4}, \frac{(2 m n)^2 + n^4 - 4 m n \sqrt{8 m^4 + n^4}}{ 16 m^4 + n^4}\right)$
\If{$x_1 \neq 1$} $ Fermat(x_1) $ \EndIf
\If{$x_2 \neq 1$} $ Fermat(x_2) $ \EndIf
\end{algorithmic}
\textbf{end procedure}
\end{algorithm} 
\newline
A python version of this algorithm can be found in the Sagemath notebook $\texttt{11\_Fermat.ipynb}$.
\pagebreak
\section{The Prime footprint tables} 
\label{section:TFP}
For the footprint equations $\ref{section:FP}$ we include the solutions $N,m,n$. 
\subsection{For $N = p ( 1 mod 8)$}
\label{table:T0}
\begin{center}
\begin{tabu}{|c|c|c|c|}  \hline
\multicolumn{2}{|c|[1pt]}{$N$}  & \multicolumn{2}{c|[1pt]}{$m,n$} \\ \hline \tabucline[1pt]\hline
\multicolumn{2}{|c|[1pt]}{41}  & \multicolumn{2}{|c|[1pt]}{ 2,1}  \\ \hline
\multicolumn{2}{|c|[1pt]}{313}  & \multicolumn{2}{|c|[1pt]}{3,2}  \\ \hline
\multicolumn{2}{|c|[1pt]}{353}  & \multicolumn{2}{|c|[1pt]}{4,1}  \\ \hline \tabucline[1pt]\hline
\multicolumn{2}{|c|[1pt]}{137}  & \multicolumn{2}{|c|[1pt]}{ 23,51}  \\ \hline
\multicolumn{2}{|c|[1pt]}{257}  & \multicolumn{2}{|c|[1pt]}{ 9,17}  \\ \hline
\multicolumn{2}{|c|[1pt]}{457}  & \multicolumn{2}{|c|[1pt]}{ 11,23}  \\ \hline
\multicolumn{2}{|c|[1pt]}{761}  & \multicolumn{2}{|c|[1pt]}{ 31,51}  \\ \hline
\end{tabu}
\end{center}
\subsection{For $N = p ( 5 mod 8)$}
\label{table:TI}
\begin{center}
\begin{tabu}{|c|c|c|c|}  \hline
\multicolumn{2}{|c|[1pt]}{$N$}  & \multicolumn{2}{c|[1pt]}{$m,n$} \\ \hline \tabucline[1pt]\hline
\multicolumn{2}{|c|[1pt]}{5}  & \multicolumn{2}{|c|[1pt]}{ 1,1}  \\ \hline
\multicolumn{2}{|c|[1pt]}{13}  & \multicolumn{2}{|c|[1pt]}{ 1,-5}  \\ \hline
\multicolumn{2}{|c|[1pt]}{29}  & \multicolumn{2}{|c|[1pt]}{ 1,-13}  \\ \hline
\multicolumn{2}{|c|[1pt]}{37}  & \multicolumn{2}{|c|[1pt]}{ 5,29}  \\ \hline
\multicolumn{2}{|c|[1pt]}{53}  & \multicolumn{2}{|c|[1pt]}{ 41,145}  \\ \hline
\multicolumn{2}{|c|[1pt]}{61}  & \multicolumn{2}{|c|[1pt]}{ 5,-89}  \\ \hline
\multicolumn{2}{|c|[1pt]}{101}  & \multicolumn{2}{|c|[1pt]}{ 53,397}  \\ \hline
\multicolumn{2}{|c|[1pt]}{109}  & \multicolumn{2}{|c|[1pt]}{ 1,5}  \\ \hline
\multicolumn{2}{|c|[1pt]}{149}  & \multicolumn{2}{|c|[1pt]}{ 1,-25}  \\ \hline
\multicolumn{2}{|c|[1pt]}{157}  & \multicolumn{2}{|c|[1pt]}{ 87005,610961}  \\ \hline
\multicolumn{2}{|c|[1pt]}{173}  & \multicolumn{2}{|c|[1pt]}{ 10865,-343141}  \\ \hline
\multicolumn{2}{|c|[1pt]}{181}  & \multicolumn{2}{|c|[1pt]}{ 13,-317}  \\ \hline
\multicolumn{2}{|c|[1pt]}{197}  & \multicolumn{2}{|c|[1pt]}{ 533,5965}  \\ \hline
\multicolumn{2}{|c|[1pt]}{229}  & \multicolumn{2}{|c|[1pt]}{ 61,-1013}  \\ \hline
\multicolumn{2}{|c|[1pt]}{269}  & \multicolumn{2}{|c|[1pt]}{ 74425,1150909}  \\ \hline
\multicolumn{2}{|c|[1pt]}{277}  & \multicolumn{2}{|c|[1pt]}{ 4264945,54532889}  \\ \hline
\multicolumn{2}{|c|[1pt]}{293}  & \multicolumn{2}{|c|[1pt]}{ 75665,-3059809}  \\ \hline
\multicolumn{2}{|c|[1pt]}{317}  & \multicolumn{2}{|c|[1pt]}{ 15073,195805}  \\ \hline
\multicolumn{2}{|c|[1pt]}{349}  & \multicolumn{2}{|c|[1pt]}{ 29,257}  \\ \hline
\multicolumn{2}{|c|[1pt]}{373}  & \multicolumn{2}{|c|[1pt]}{ 312842465,-8322284921}  \\ \hline
\multicolumn{2}{|c|[1pt]}{389}  & \multicolumn{2}{|c|[1pt]}{ 976825,15047929}  \\ \hline
\multicolumn{2}{|c|[1pt]}{397}  & \multicolumn{2}{|c|[1pt]}{ 113,965}  \\ \hline
\multicolumn{2}{|c|[1pt]}{421}  & \multicolumn{2}{|c|[1pt]}{ 1453,22565}  \\ \hline
\multicolumn{2}{|c|[1pt]}{461}  & \multicolumn{2}{|c|[1pt]}{ 84653,-3918457}  \\ \hline
\multicolumn{2}{|c|[1pt]}{509}  & \multicolumn{2}{|c|[1pt]}{ 1,13}  \\ \hline
\multicolumn{2}{|c|[1pt]}{541}  & \multicolumn{2}{|c|[1pt]}{ 40573,-2037137}  \\ \hline
\multicolumn{2}{|c|[1pt]}{557}  & \multicolumn{2}{|c|[1pt]}{ 5585,-312709}  \\ \hline
\multicolumn{2}{|c|[1pt]}{613}  & \multicolumn{2}{|c|[1pt]}{ 823709,17966645}  \\ \hline
\multicolumn{2}{|c|[1pt]}{653}  & \multicolumn{2}{|c|[1pt]}{ 669773012805221,7399743723764065}  \\ \hline
\multicolumn{2}{|c|[1pt]}{661}  & \multicolumn{2}{|c|[1pt]}{ 923777,-56911865}  \\ \hline
\multicolumn{2}{|c|[1pt]}{677}  & \multicolumn{2}{|c|[1pt]}{ 4755629897398453,55948826401549645}  \\ \hline
\multicolumn{2}{|c|[1pt]}{701}  & \multicolumn{2}{|c|[1pt]}{ 265358113,2932361965}  \\ \hline
\multicolumn{2}{|c|[1pt]}{709}  & \multicolumn{2}{|c|[1pt]}{ 1,17}  \\ \hline
\multicolumn{2}{|c|[1pt]}{733}  & \multicolumn{2}{|c|[1pt]}{ 24668160920485,307594098561689}  \\ \hline
\multicolumn{2}{|c|[1pt]}{757}  & \multicolumn{2}{|c|[1pt]}{ 1607466650773733245,-83390357142888376589}  \\ \hline
\multicolumn{2}{|c|[1pt]}{773}  & \multicolumn{2}{|c|[1pt]}{ 1528361,23404585}  \\ \hline
\multicolumn{2}{|c|[1pt]}{797}  & \multicolumn{2}{|c|[1pt]}{ 18638233339062102833,349269016577216709805}  \\ \hline
\multicolumn{2}{|c|[1pt]}{821}  & \multicolumn{2}{|c|[1pt]}{ 5,-181}  \\ \hline
\multicolumn{2}{|c|[1pt]}{829}  & \multicolumn{2}{|c|[1pt]}{ 240841442845,4628302886849}  \\ \hline
\multicolumn{2}{|c|[1pt]}{853}  & \multicolumn{2}{|c|[1pt]}{ 4849409,-320740745}  \\ \hline
\multicolumn{2}{|c|[1pt]}{877}  & \multicolumn{2}{|c|[1pt]}{ 54778385345,-1789378941029}  \\ \hline
\multicolumn{2}{|c|[1pt]}{941}  & \multicolumn{2}{|c|[1pt]}{ 738533,9409777}  \\ \hline
\multicolumn{2}{|c|[1pt]}{997}  & \multicolumn{2}{|c|[1pt]}{ 93521216025186467821517,2200234831458939782292245}  \\ \hline  
\end{tabu}
\end{center}
\pagebreak
\subsection{For $N = p ( 7 mod 8)$ }
\label{table:TII}
\begin{center}
\begin{tabu}{|c|c|c|c|}  \hline
\multicolumn{2}{|c|[1pt]}{$N$}  & \multicolumn{2}{|c|[1pt]}{$m,n$} \\ \hline \tabucline[1pt]\hline
\multicolumn{2}{|c|[1pt]}{7}  & \multicolumn{2}{|c|[1pt]}{ 2,1}  \\ \hline
\multicolumn{2}{|c|[1pt]}{23}  & \multicolumn{2}{|c|[1pt]}{ 13,6}  \\ \hline
\multicolumn{2}{|c|[1pt]}{31}  & \multicolumn{2}{|c|[1pt]}{ 5,4}  \\ \hline
\multicolumn{2}{|c|[1pt]}{47}  & \multicolumn{2}{|c|[1pt]}{ 53,36}  \\ \hline
\multicolumn{2}{|c|[1pt]}{71}  & \multicolumn{2}{|c|[1pt]}{ 6,5}  \\ \hline
\multicolumn{2}{|c|[1pt]}{79}  & \multicolumn{2}{|c|[1pt]}{ 125,52}  \\ \hline
\multicolumn{2}{|c|[1pt]}{103}  & \multicolumn{2}{|c|[1pt]}{ 10361,4522}  \\ \hline
\multicolumn{2}{|c|[1pt]}{127}  & \multicolumn{2}{|c|[1pt]}{ 145249,60248}  \\ \hline
\multicolumn{2}{|c|[1pt]}{151}  & \multicolumn{2}{|c|[1pt]}{ 17,10}  \\ \hline
\multicolumn{2}{|c|[1pt]}{167}  & \multicolumn{2}{|c|[1pt]}{ 449,378}  \\ \hline
\multicolumn{2}{|c|[1pt]}{191}  & \multicolumn{2}{|c|[1pt]}{ 27469,11580}  \\ \hline
\multicolumn{2}{|c|[1pt]}{199}  & \multicolumn{2}{|c|[1pt]}{ 1465,1274}  \\ \hline
\multicolumn{2}{|c|[1pt]}{223}  & \multicolumn{2}{|c|[1pt]}{ 77056,46313}  \\ \hline
\multicolumn{2}{|c|[1pt]}{239}  & \multicolumn{2}{|c|[1pt]}{ 12,5}  \\ \hline
\multicolumn{2}{|c|[1pt]}{263}  & \multicolumn{2}{|c|[1pt]}{ 239558789,102570078}  \\ \hline
\multicolumn{2}{|c|[1pt]}{271}  & \multicolumn{2}{|c|[1pt]}{ 14485,10388}  \\ \hline
\multicolumn{2}{|c|[1pt]}{311}  & \multicolumn{2}{|c|[1pt]}{ 676273,322890}  \\ \hline
\multicolumn{2}{|c|[1pt]}{359}  & \multicolumn{2}{|c|[1pt]}{ 8329,3450}  \\ \hline
\multicolumn{2}{|c|[1pt]}{367}  & \multicolumn{2}{|c|[1pt]}{ 91887589669,89142445028}  \\ \hline
\multicolumn{2}{|c|[1pt]}{383}  & \multicolumn{2}{|c|[1pt]}{ 4692,2141}  \\ \hline
\multicolumn{2}{|c|[1pt]}{431}  & \multicolumn{2}{|c|[1pt]}{ 15805,13908}  \\ \hline
\multicolumn{2}{|c|[1pt]}{439}  & \multicolumn{2}{|c|[1pt]}{ 16510,8069}  \\ \hline
\multicolumn{2}{|c|[1pt]}{463}  & \multicolumn{2}{|c|[1pt]}{ 27557,19564}  \\ \hline
\multicolumn{2}{|c|[1pt]}{479}  & \multicolumn{2}{|c|[1pt]}{ 193,120}  \\ \hline
\multicolumn{2}{|c|[1pt]}{487}  & \multicolumn{2}{|c|[1pt]}{ 5386013,2963446}  \\ \hline
\multicolumn{2}{|c|[1pt]}{503}  & \multicolumn{2}{|c|[1pt]}{ 388326872921,162262911162}  \\ \hline
\multicolumn{2}{|c|[1pt]}{599}  & \multicolumn{2}{|c|[1pt]}{ 1317481516325,661578829566}  \\ \hline
\multicolumn{2}{|c|[1pt]}{607}  & \multicolumn{2}{|c|[1pt]}{ 94040069519467036,59514405263901653}  \\ \hline
\multicolumn{2}{|c|[1pt]}{631}  & \multicolumn{2}{|c|[1pt]}{ 10078210,4186673}  \\ \hline
\multicolumn{2}{|c|[1pt]}{647}  & \multicolumn{2}{|c|[1pt]}{ 555349528702900386,1035251108850383833}  \\ \hline
\multicolumn{2}{|c|[1pt]}{719}  & \multicolumn{2}{|c|[1pt]}{ 65,144}  \\ \hline
\multicolumn{2}{|c|[1pt]}{727}  & \multicolumn{2}{|c|[1pt]}{ 3971654742468682789,2041107876866928758}  \\ \hline
\multicolumn{2}{|c|[1pt]}{743}  & \multicolumn{2}{|c|[1pt]}{ 607763379942529,546999398401098}  \\ \hline
\multicolumn{2}{|c|[1pt]}{751}  & \multicolumn{2}{|c|[1pt]}{ 20,13}  \\ \hline
\multicolumn{2}{|c|[1pt]}{823}  & \multicolumn{2}{|c|[1pt]}{ 39999941080334034526,22509182182068729173}  \\ \hline
\multicolumn{2}{|c|[1pt]}{839}  & \multicolumn{2}{|c|[1pt]}{ 1050,457}  \\ \hline
\multicolumn{2}{|c|[1pt]}{863}  & \multicolumn{2}{|c|[1pt]}{ 40579118388,16831312549}  \\ \hline
\multicolumn{2}{|c|[1pt]}{887}  & \multicolumn{2}{|c|[1pt]}{ 37615957468559969,28312289691904218}  \\ \hline
\multicolumn{2}{|c|[1pt]}{911}  & \multicolumn{2}{|c|[1pt]}{ 349810551819545,145327242840624}  \\ \hline
\multicolumn{2}{|c|[1pt]}{919}  & \multicolumn{2}{|c|[1pt]}{ 20306,10225}  \\ \hline
\multicolumn{2}{|c|[1pt]}{967}  & \multicolumn{2}{|c|[1pt]}{ 6445185278237,3169446617854}  \\ \hline
\multicolumn{2}{|c|[1pt]}{983}  & \multicolumn{2}{|c|[1pt]}{ 14410159497869814,6444783543384757}  \\ \hline
\multicolumn{2}{|c|[1pt]}{991}  & \multicolumn{2}{|c|[1pt]}{ 659581,282740}  \\ \hline
\end{tabu}
\end{center}
\begin{center}  
\renewcommand{\arraystretch}{1}
\end{center}
\pagebreak
\subsection{For $N = 2  p , p (7 mod 8)$}
\label{table:TIII}
\begin{center}
\begin{tabu}{|c|c|c|c|}  \hline
\multicolumn{2}{|c|[1pt]}{$N$}  & \multicolumn{2}{|c|[1pt]}{$m,n$} \\ \hline \tabucline[1pt]\hline
\multicolumn{2}{|c|[1pt]}{14}   & \multicolumn{2}{|c|[1pt]}{$ 2,1$} \\ \hline
\multicolumn{2}{|c|[1pt]}{46}   & \multicolumn{2}{|c|[1pt]}{$ 2,-3$} \\ \hline
\multicolumn{2}{|c|[1pt]}{62}  & \multicolumn{2}{|c|[1pt]}{$ 20,7$} \\ \hline
\multicolumn{2}{|c|[1pt]}{94}  & \multicolumn{2}{|c|[1pt]}{$ 12,7$} \\ \hline
\multicolumn{2}{|c|[1pt]}{142}  & \multicolumn{2}{|c|[1pt]}{$ 290,171$} \\ \hline
\multicolumn{2}{|c|[1pt]}{158}  & \multicolumn{2}{|c|[1pt]}{$ 20,-31$} \\ \hline
\multicolumn{2}{|c|[1pt]}{206}  & \multicolumn{2}{|c|[1pt]}{$ 34,-41$} \\ \hline
\multicolumn{2}{|c|[1pt]}{254}  & \multicolumn{2}{|c|[1pt]}{$ 16,-17$} \\ \hline
\multicolumn{2}{|c|[1pt]}{302}  & \multicolumn{2}{|c|[1pt]}{$ 35890,-60401$} \\ \hline
\multicolumn{2}{|c|[1pt]}{334}  & \multicolumn{2}{|c|[1pt]}{$ 1018,-1257$} \\ \hline
\multicolumn{2}{|c|[1pt]}{382}  & \multicolumn{2}{|c|[1pt]}{$ 540,239$} \\ \hline
\multicolumn{2}{|c|[1pt]}{398}  & \multicolumn{2}{|c|[1pt]}{$ 650,-1009$} \\ \hline
\multicolumn{2}{|c|[1pt]}{446}  & \multicolumn{2}{|c|[1pt]}{$ 6104,-6169$} \\ \hline
\multicolumn{2}{|c|[1pt]}{478}  & \multicolumn{2}{|c|[1pt]}{$ 3481500,2186791$} \\ \hline
\multicolumn{2}{|c|[1pt]}{526}  & \multicolumn{2}{|c|[1pt]}{$ 74,51$} \\ \hline
\multicolumn{2}{|c|[1pt]}{542}  & \multicolumn{2}{|c|[1pt]}{$ 4397260,1896071$} \\ \hline
\multicolumn{2}{|c|[1pt]}{622}  & \multicolumn{2}{|c|[1pt]}{$ 10,-9$} \\ \hline
\multicolumn{2}{|c|[1pt]}{718}  & \multicolumn{2}{|c|[1pt]}{$ 65870050,-68694393$} \\ \hline
\multicolumn{2}{|c|[1pt]}{734}  & \multicolumn{2}{|c|[1pt]}{$ 28,-23$} \\ \hline
\multicolumn{2}{|c|[1pt]}{766}  & \multicolumn{2}{|c|[1pt]}{$ 26316,-22351$} \\ \hline
\multicolumn{2}{|c|[1pt]}{862}  & \multicolumn{2}{|c|[1pt]}{$ 10373408340,3061104191$} \\ \hline
\multicolumn{2}{|c|[1pt]}{878}  & \multicolumn{2}{|c|[1pt]}{$ 2290,683$} \\ \hline
\multicolumn{2}{|c|[1pt]}{926}  & \multicolumn{2}{|c|[1pt]}{$ 19163084,10851559$} \\ \hline
\multicolumn{2}{|c|[1pt]}{958}  & \multicolumn{2}{|c|[1pt]}{$ 196826640,-325857553$} \\ \hline
\multicolumn{2}{|c|[1pt]}{974}  & \multicolumn{2}{|c|[1pt]}{$ 1451962,833243$} \\ \hline
\end{tabu}
\end{center}
\subsection{For $N = 2  p , p (3 mod 8)$}
\label{table:TIV}
\begin{center}
\begin{tabu}{|c|c|c|c|}  \hline
\multicolumn{2}{|c|[1pt]}{$N$}  & \multicolumn{2}{|c|[1pt]}{$m,n$} \\ \hline \tabucline[1pt]\hline
\multicolumn{2}{|c|[1pt]}{6}  & \multicolumn{2}{|c|[1pt]}{$ 1,0$} \\ \hline
\multicolumn{2}{|c|[1pt]}{22}  & \multicolumn{2}{|c|[1pt]}{$ 2,1$} \\ \hline
\multicolumn{2}{|c|[1pt]}{38}  & \multicolumn{2}{|c|[1pt]}{$ 3,4$} \\ \hline
\multicolumn{2}{|c|[1pt]}{86}  & \multicolumn{2}{|c|[1pt]}{$ 2,3$} \\ \hline
\multicolumn{2}{|c|[1pt]}{118}  & \multicolumn{2}{|c|[1pt]}{$ 131,-42$} \\ \hline
\multicolumn{2}{|c|[1pt]}{134}  & \multicolumn{2}{|c|[1pt]}{$ 30,-11$} \\ \hline
\multicolumn{2}{|c|[1pt]}{166}  & \multicolumn{2}{|c|[1pt]}{$ 70,71$} \\ \hline
\multicolumn{2}{|c|[1pt]}{214}  & \multicolumn{2}{|c|[1pt]}{$ 12,13$} \\ \hline
\multicolumn{2}{|c|[1pt]}{262}  & \multicolumn{2}{|c|[1pt]}{$ 71,48$} \\ \hline
\multicolumn{2}{|c|[1pt]}{278}  & \multicolumn{2}{|c|[1pt]}{$ 767,336$} \\ \hline
\multicolumn{2}{|c|[1pt]}{326}  & \multicolumn{2}{|c|[1pt]}{$ 170,-69$} \\ \hline
\multicolumn{2}{|c|[1pt]}{358}  & \multicolumn{2}{|c|[1pt]}{$ 4674,-143$} \\ \hline
\multicolumn{2}{|c|[1pt]}{422}  & \multicolumn{2}{|c|[1pt]}{$ 23,-6$} \\ \hline
\multicolumn{2}{|c|[1pt]}{454}  & \multicolumn{2}{|c|[1pt]}{$ 700,467$} \\ \hline
\multicolumn{2}{|c|[1pt]}{502}  & \multicolumn{2}{|c|[1pt]}{$ 3529,-848$} \\ \hline
\multicolumn{2}{|c|[1pt]}{566}  & \multicolumn{2}{|c|[1pt]}{$ 769,1590$} \\ \hline
\multicolumn{2}{|c|[1pt]}{614}  & \multicolumn{2}{|c|[1pt]}{$ 7741,-174$} \\ \hline
\multicolumn{2}{|c|[1pt]}{662}  & \multicolumn{2}{|c|[1pt]}{$ 167376,128473$} \\ \hline
\multicolumn{2}{|c|[1pt]}{694}  & \multicolumn{2}{|c|[1pt]}{$ 13,10$} \\ \hline
\multicolumn{2}{|c|[1pt]}{758}  & \multicolumn{2}{|c|[1pt]}{$ 146000702,316270911$} \\ \hline
\multicolumn{2}{|c|[1pt]}{838}  & \multicolumn{2}{|c|[1pt]}{$ 2690159,-878878$} \\ \hline
\multicolumn{2}{|c|[1pt]}{886}  & \multicolumn{2}{|c|[1pt]}{$ 1334,181$} \\ \hline
\multicolumn{2}{|c|[1pt]}{934}  & \multicolumn{2}{|c|[1pt]}{$ 11,-4$} \\ \hline
\multicolumn{2}{|c|[1pt]}{982}  & \multicolumn{2}{|c|[1pt]}{$ 56212854,-14198593$} \\ \hline
\multicolumn{2}{|c|[1pt]}{998}  & \multicolumn{2}{|c|[1pt]}{$ 1839056773,1911574194$} \\ \hline
\end{tabu}
\end{center}
\pagebreak
\section{Conclusions}
The fascination for Fermat's solution and the Pythagorean theorem led to the rational triples their connection to the congruent number problem and their discoveries going beyond like the Trinity system unfolding beautiful algebraic and geometric properties.
\newline
\includegraphics[width=150mm,scale=1]{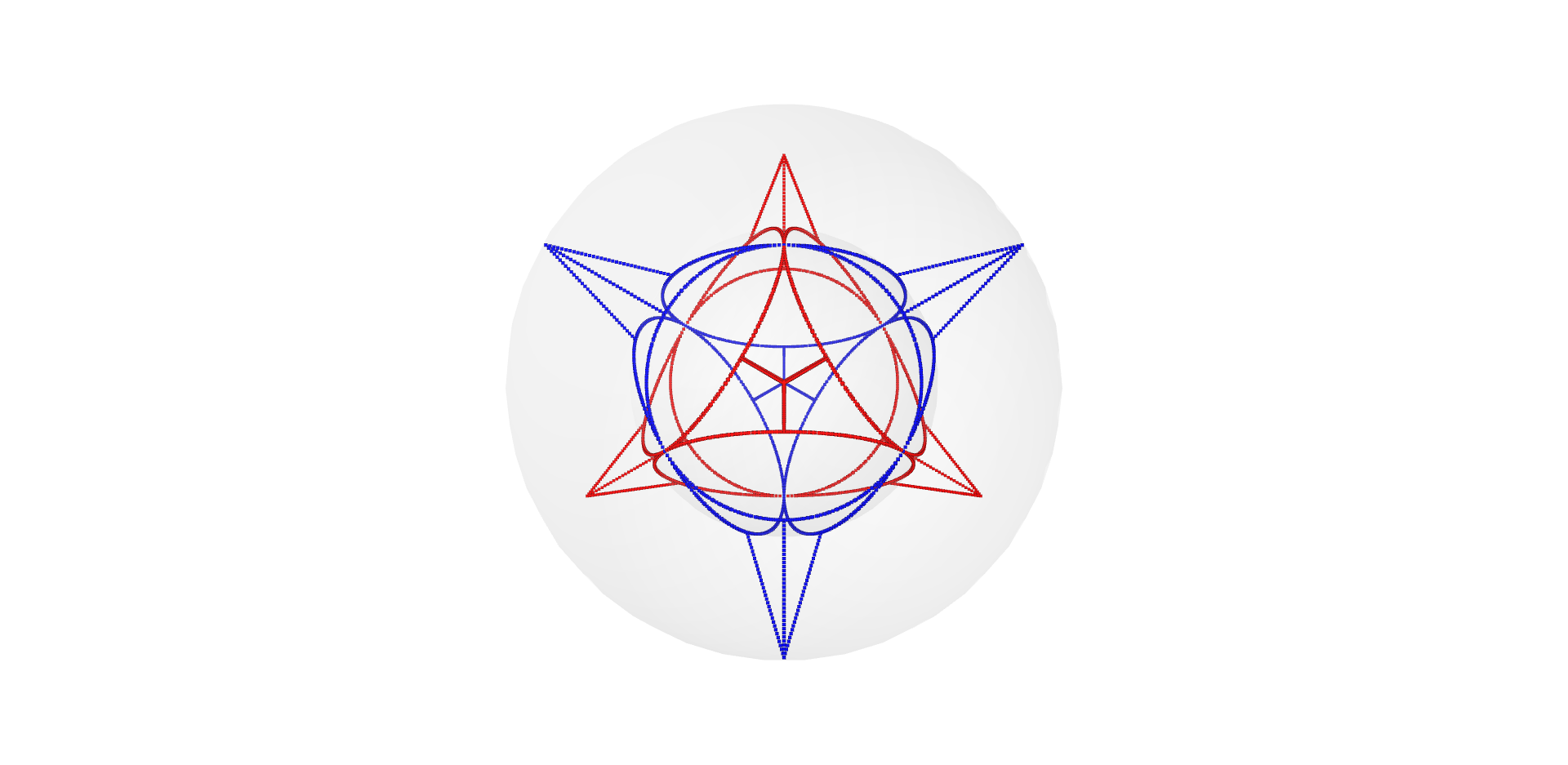}
\section{Acknowledgment}
I would like to express my deepest gratitude to Dr. N.Tati Ruiz for supporting and motivating me.
\bibliographystyle{te}
\bibliography{bibl}

\end{document}